\documentclass[11pt]{amsart}
\usepackage{amsmath,amssymb,amsthm,amsrefs,geometry,graphics,color,mdwlist,mathrsfs,enumitem,mathtools,tikz-cd,amscd,mathrsfs,tikz,tabstackengine,nicematrix,graphicx,xcolor,varwidth}
\usepackage[all]{xy}

\usepackage{hyperref}
\hypersetup{hypertex=true,
colorlinks=true,
linkcolor=blue,
anchorcolor=blue,
citecolor=blue}

\numberwithin{equation}{section}
\title{Singularity Categories of B{\"a}ckstr{\"o}m Orders}
\author{HONGRUI WEI}

\address{Graduate School of Mathematical Sciences, The University of Tokyo, Komaba 3-8-1, Meguro-ku, Tokyo, 153-8941, Japan}
\email{hongrui-wei@g.ecc.u-tokyo.ac.jp}

\begin{document}
\theoremstyle{plain}
\newtheorem{theorem}{Theorem}[section]
\newtheorem{proposition}[theorem]{Proposition}
\newtheorem{lemma}[theorem]{Lemma}
\newtheorem{conjecture}[theorem]{Conjecture}
\newtheorem{corollary}[theorem]{Corollary}
\newtheorem{axiom}[theorem]{Axiom}
\newtheorem{thmdef}[theorem]{Theorem-Definition}
\newtheorem{defprop}[theorem]{Definition-Proposition}
\newtheorem{slogan}[theorem]{Slogan}

\theoremstyle{definition}
\newtheorem{definition}[theorem]{Definition}
\newtheorem{example}[theorem]{Example}
\newtheorem*{assumption}{Assumption}
\newtheorem*{assumptions}{Assumptions}
\newtheorem{notation}[theorem]{Notation}
\newtheorem{question}[theorem]{Question}
\newtheorem*{conventions and notations}{Conventions and Notations}
\newtheorem*{acknowledgement}{Acknowledgement}
\newtheorem*{claim}{Claim}

\theoremstyle{remark}
\newtheorem{remark}[theorem]{Remark}
\newtheorem{summary}[theorem]{Summary}
\newtheorem{observation}[theorem]{Observation}
\newtheorem{conclusion}[theorem]{Conclusion}
\newtheorem{case}[theorem]{Case}

\newcommand{\Hom}{\operatorname{Hom}}
\newcommand{\sHom}{\operatorname{\underline{Hom}}}
\newcommand{\cHom}{\operatorname{\overline{Hom}}}
\newcommand{\End}{\operatorname{End}}
\newcommand{\sEnd}{\operatorname{\underline{End}}}
\newcommand{\Ext}{\operatorname{Ext}}
\newcommand{\Tor}{\operatorname{Tor}}

\renewcommand{\Im}{\operatorname{Im}}
\newcommand{\Ker}{\operatorname{Ker}}
\newcommand{\Coker}{\operatorname{Coker}}
\newcommand{\Id}{\operatorname{Id}}
\newcommand{\res}{\operatorname{res}}

\newcommand{\pd}{\operatorname{pd}}
\newcommand{\id}{\operatorname{id}}
\newcommand{\gd}{\operatorname{gl.dim}}
\newcommand{\Gpd}{\operatorname{Gpd}}

\newcommand{\obj}{\operatorname{Obj}}
\newcommand{\op}{\operatorname{op}}
\newcommand{\ind}{\operatorname{ind}}
\newcommand{\simple}{\operatorname{sim}}
\renewcommand{\mod}{\operatorname{mod}}
\newcommand{\smod}{\operatorname{\underline{mod}}}
\newcommand{\cmod}{\operatorname{\overline{mod}}}
\newcommand{\Frac}{\operatorname{Frac}}
\newcommand{\rad}{\operatorname{rad}}
\newcommand{\corad}{\operatorname{corad}}
\newcommand{\soc}{\operatorname{soc}}
\newcommand{\fl}{\operatorname{f.l.}}
\newcommand{\CM}{\operatorname{CM}}
\newcommand{\sCM}{\operatorname{\underline{CM}}}
\newcommand{\Dsg}{\operatorname{\rm{D}_{sg}}}
\newcommand{\proj}{\operatorname{proj}}
\newcommand{\Gp}{\operatorname{Gproj}}
\newcommand{\sGp}{\operatorname{\underline{Gproj}}}
\newcommand{\inj}{\operatorname{inj}}
\newcommand{\add}{\operatorname{add}}
\newcommand{\sadd}{\operatorname{\underline{add}}}

\newcommand{\depth}{\operatorname{depth}}
\renewcommand{\dim}{\operatorname{dim}}

\newcommand{\Gen}{\operatorname{Gen}}
\newcommand{\Cogen}{\operatorname{Cogen}}

\newcommand{\Spec}{\operatorname{Spec}}
\newcommand{\Tr}{\operatorname{Tr}}

\newcommand{\N}{\mathbb{N}}
\newcommand{\Z}{\mathbb{Z}}
\newcommand{\Q}{\mathbb{Q}}
\newcommand{\R}{\mathbb{R}}
\newcommand{\C}{\mathbb{C}}

\newcommand{\Sp}{\operatorname{Sp}}
\newcommand{\rep}{\operatorname{rep}}

\newcommand{\rightarrowdbl}{\rightarrow\mathrel{\mkern-14mu}\rightarrow}
\newcommand{\xrightarrowdbl}[2][]{%
  \xrightarrow[#1]{#2}\mathrel{\mkern-14mu}\rightarrow
}


\begin{abstract}
    B{\"a}ckstr{\"o}m orders are a class of algebras over complete discrete valuation rings. Their Cohen-Macaulay representations are in correspondence with the representations of certain quivers/species by Ringel and Roggenkamp. In this paper, we give explicit descriptions of the singularity categories of B{\"a}ckstr{\"o}m orders via certain von Neumann regular algebras and associated bimodules. We further provide singular equivalences between B{\"a}ckstr{\"o}m orders and specific finite dimensional radical square zero algebras. We also classify weakly regular, Gorenstein, Iwanaga-Gorenstein and sg-Hom-finite B{\"a}ckstr{\"o}m orders.
\end{abstract}

\maketitle
\setcounter{tocdepth}{2}

\tableofcontents

\section{Introduction}

\subsection{Background} \label{sec:motivation}

\subsubsection{Cohen-Macaulay representations} 
The study of Cohen-Macaulay modules is an active branch of representation theory of associative algebras, having an important connection with algebraic geometry and physics \cites{Yos90, Sim92,LW12, IW14}. By definition, an \textit{$R$-order} $\Lambda$ is an $R$-algebra over a complete regular local ring which is free of finite rank over $R$. This notion generalizes finite dimensional algebras to the case where the base rings have positive Krull dimension. Cohen-Macaulay modules are widely studied for such class of algebras, especially when $\dim R=1$ \cites{CR90, HN94, Rei03}.

\subsubsection{Semisimple algebras and hereditary orders} 
A basic class of orders is called \textit{regular orders} (see Figure (\ref{eqa:hierarchy})), i.e. $\gd \Lambda=\dim R$. For instance, if $R$ is a field, then $\Lambda$ is regular if and only if $\Lambda$ is semisimple; if $R$ is a complete discrete valuation ring, then $\Lambda$ is regular if and only if $\Lambda$ is hereditary. In this case, $\CM \Lambda = \proj \Lambda$. The well-known Artin-Wedderburn theorem characterizes regular orders over fields. A similar characterization holds for regular orders over complete discrete valuation rings (see Theorem \ref{thm:structure theorem of hereditary orders}).

\subsubsection{Radical square zero algebras and B{\"a}ckstr{\"o}m orders}
Semisimple algebras are also characterized by the property that their Jacobson radicals are trivial. From this perspective, finite dimensional radical square zero algebras are close to semisimple algebras. We have $\ind (\mod A) \xleftrightarrow{1:1} \ind (\mod H')$ for a radical square zero algebra $A$,  where 
\begin{equation*}
    H'=
    \begin{pmatrix}
        A / \rad A & 0 \\
        \rad A & A / \rad A
    \end{pmatrix}.
\end{equation*} 
The representation theory of $A$ is tractable since $H$ is a finite dimensional hereditary algebra, whose representations correspond to the representations of a certain acyclic quiver/species. Moreover, a stable equivalence exists: $\smod A \simeq \smod H'$. See \cite{Gab73} and \cite[Section X.2]{ARS97}.  

For $\dim R=1$, a class of orders which are ideal-theoretically close to hereditary orders is called B{\"a}ckstr{\"o}m orders. By definition, a \textit{B{\"a}ckstr{\"o}m order} $\Lambda=(\Lambda,\Gamma)$ contains a pair of orders $\Lambda\subseteq \Gamma$ such that $\Gamma$ is a hereditary order and $\rad \Lambda =\rad \Gamma$; see Definition \ref{def:Backorder}. They include several notable families of orders arising naturally in topology and geometry, such as ribbon graph orders \cites{KR01,Gne19}, gentle orders \cite{Gne25}, nodal orders \cites{BD24a,BD24b} and quadratic orders \cite{Iya05}. Furthermore, this framework has been fruitfully extended to the broader notion of Bäckström rings \cite{Dro23}, which contain important classes of algebras such as gentle and skew-gentle algebras.

The representation theory of B{\"a}ckstr{\"o}m orders was developed by Ringel and Roggenkamp \cite{RR79} (see also Green and Reiner \cite{GR78}). We have $\ind (\CM \Lambda) \xleftrightarrow{\text{almost 1:1}} \ind (\mod H)$, where 
\begin{equation*}
    H=
    \begin{pmatrix}
        \Gamma / \rad \Gamma & 0 \\
        \Gamma / \rad \Gamma & \Lambda / \rad \Lambda 
    \end{pmatrix}
\end{equation*}
is a finite dimensional hereditary algebra; see Theorem \ref{thm:rep_of_Back}. Based on this correspondence, we provide a categorical equivalence. 
\begin{theorem}[=Theorem \ref{thm:stable equivalence}] \label{thmintro:stable equivalence}
    Assume that $\Lambda=(\Lambda,\Gamma)$ is a B{\"a}ckstr{\"o}m order and keep the notations above. Then we have a categorical equivalence $\sCM \Lambda \simeq \smod_s H$, where $\smod_s H$ is the full subcategory of $\smod H$ whose objects do not contain any non-zero simple direct summand.
\end{theorem}

We have seen that the representation theory of radical square zero algebras and B{\"a}ckstr{\"o}m orders behaves similarly: both correspond to hereditary algebras of the form of upper triangle matrices. Motivated by this phenomenon, we purpose the following slogan:

\begin{slogan}\label{slogan}
    B{\"a}ckstr{\"o}m orders are one-dimensional counterparts of finite dimensional radical square zero algebras.
\end{slogan}

In this paper, we focus on the singularity category of B{\"a}ckstr{\"o}m orders and provide further evidence to support this slogan.

\subsection{Main Results}

\subsubsection{Singularity categories of B{\"a}ckstr{\"o}m orders and radical square zero algebras}
The \textit{singularity category} $\Dsg(\Lambda)$ of a noetherian ring $\Lambda$ was introduced by Buchweitz \cite{Buc97}. It is defined as the Verdier quotient of the bounded derived category $\rm{D}^{\rm{b}}(\Lambda)$ of $\mod \Lambda$ by the bounded homotopy category $\rm{K}^{\rm{b}}(\proj \Lambda)$. The singularity category captures the stable homological property of a noetherian ring. More precisely, as shown in \cite{KV87}, the singularity category is equivalent to the \textit{stabilization} $\mathcal{S}(\smod\Lambda, \Omega)$ of the stable module category $\smod \Lambda$ with the syzygy functor $\Omega$ on it, where the morphism sets are direct limits of the syzygy functor; see Section \ref{subsec: Stabilizations and Singularity Categories}. We will use this fact to study and give an explicit description of the singularity category of a B{\"a}ckstr{\"o}m order. To illustrate our result, we introduce our settings.

Assume that $\Lambda=(\Lambda,\Gamma)$ is a B{\"a}ckstr{\"o}m order. Let $\Ker \mu$ be the kernel of the multiplication map $\mu:\Gamma \otimes_R \Lambda \rightarrowdbl \Gamma$. Let $D:=\sEnd_\Lambda (\Gamma)$ and $M:=\sHom_\Lambda(\Gamma,\Ker \mu)$. Then $M$ has a $(D,D)$-bimodule structure which induces a functor $?\otimes_D M: \mod D \to \mod D$; see Proposition \ref{prop: M is a (D,D)-bimodule}. Let $M^{\otimes 0}:=D$ and $M^{\otimes i}:=M\otimes_D M^{\otimes i-1}$ for $i\geq1$. Define 
    \begin{equation*}
        V(\Lambda):=\underset{i\geq0}{\varinjlim} \End_D(M^{\otimes i}) \quad {\rm{and}} \quad 
        K(\Lambda):=\underset{i\geq1}{\varinjlim} \Hom_D(M^{\otimes i}, M^{\otimes i-1}).
    \end{equation*}
The following is our first main result.

\begin{theorem}[=Theorem \ref{thm:DsgBack}]\label{thm_intro:DsgBack}
    Assume that $\Lambda=(\Lambda,\Gamma)$ is a B{\"a}ckstr{\"o}m order. 
    Then $V(\Lambda)$ is a von Neumann regular algebra (see \rm{{Theorem-Definition} \ref{thmdef:von Neumann regular}}) and there are triangle equivalences
    \begin{equation*}
        (\Dsg(\Lambda), [1]) \simeq \mathcal{S}(\smod \Lambda, \Omega) \simeq (\proj V(\Lambda),? \otimes_{V(\Lambda)} K(\Lambda)),
    \end{equation*}
    where the triangulated structures are trivial (see \rm{Section \ref{subsec: Stabilizations and Singularity Categories}}).
\end{theorem}

For a finite dimensional radical square zero algebra, a similar description was obtained by Xiao-Wu Chen using the theory of stabilization \cite{Che11}; see Theorem \ref{thm:DsgArtin}. We can see that our result and Chen's result exhibit a close resemblance to each other, which provides another evidence for the guiding analogy in our slogan.

\subsubsection{Singular equivalences between B{\"a}ckstr{\"o}m orders and radical square zero algebras}
By a \textit{singular equivalence} between two noetherian rings, we mean a triangle equivalence between their singularity categories. It is a weaker version of derived  equivalence. Combining together Theorem \ref{thm_intro:DsgBack} and Chen's result, we construct a singular equivalence between a B{\"a}ckstr{\"o}m order and a finite dimensional radical square zero algebra. This provides another connection between these two classes of algebras.

\begin{corollary}[=Corollary \ref{cor:DsgLambda=DsgA}]\label{corintro:DsgLambda=DsgA}
    Let $\Lambda=(\Lambda,\Gamma)$ be a B{\"a}ckstr{\"o}m order and keep the notations as in {\rm{Theorem \ref{thm_intro:DsgBack}}}. Define $A(\Lambda)$ as the trivial extension 
    \begin{equation*}
        A(\Lambda)=D\oplus M,
    \end{equation*}
    which is a finite dimensional radical square zero algebra. Then we have a triangle equivalence
    \begin{equation*}
        \Dsg(\Lambda)\simeq\Dsg(A(\Lambda)).
    \end{equation*} 
\end{corollary}

\subsubsection{Classifications of some classes of B{\"a}ckstr{\"o}m orders}
Within the framework of orders, there are several notable classes whose hierarchy is presented below. Let $\Lambda$ be an $R$-order. Then we have the following hierarchy:
\begin{equation}\label{eqa:hierarchy}
    \begin{tikzcd}[column sep={between origins, 14em}, row sep=large]
        \fbox{\begin{tabular}{c} Regular\\ $(\gd \Lambda=\dim R)$ \end{tabular}}  \arrow[r,Rightarrow] \arrow[d,Rightarrow,to path={([yshift=-3.8em]\tikztostart.north) -- ([yshift=0em]\tikztotarget.north)}] & \fbox{\begin{tabular}{@{}c@{}} Gorenstein\\ $(\id _{\Lambda}\Lambda=\id_{\Lambda^{\op}} \Lambda=\dim R)$ \end{tabular}} \arrow[d,Rightarrow] & \\
        \fbox{\begin{tabular}{c} Weakly regular \\ ($\gd \Lambda < \infty)$\end{tabular}} \arrow[r,Rightarrow, to path={([xshift=0.8em]\tikztostart.east) -- ([xshift=-0.35em]\tikztotarget.west)}] & \fbox{\begin{tabular}{c} Iwanaga-Gorenstein\\ $(\id _{\Lambda}\Lambda=\id_{\Lambda^{\op}} \Lambda< \infty)$ \end{tabular}} \arrow[r,Rightarrow, to path={([xshift=0.2em]\tikztostart.east) -- ([xshift=-0.25em]\tikztotarget.west)}]& \fbox{\begin{tabular}{c} sg-Hom-finite \\ ($\Dsg(\Lambda)$ is Hom-finite) \end{tabular}} 
    \end{tikzcd}
\end{equation}

If $R$ is a field, then $\Lambda$ is Gorenstein if and only if $\Lambda$ is self-injective. The singular equivalence in Corollary \ref{corintro:DsgLambda=DsgA} additionally preserves some important homological properties. As application, we use it to classify some important classes of B{\"a}ckstr{\"o}m orders. For the finite dimensional algebra $A(\Lambda)$, we may assign a \textit{valued quiver} $Q_{A(\Lambda)}$ to it; see Section \ref{subsec: Classifications of some Classes of B-Orders}. 

\begin{theorem}[=Corollary \ref{cor:HomfiniteDsgBack}, Corollary \ref{cor:glfiniteBack}, Theorem \ref{cor:restriction to GP}, Theorem \ref{cor:criterion for Gor}] \label{thmintro: criterion for classes of B-orders}
    Let $\Lambda$ be a B{\"a}ckstr{\"o}m order, $A(\Lambda)$ the finite dimensional radical square zero algebra associated to it and $Q_{A(\Lambda)}$ the valued quiver of $A(\Lambda)$. We have the following criteria.
    \begin{enumerate}[label=\textnormal{(\arabic*)}]
        \item $\Lambda$ is weakly regular if and only if $A(\Lambda)$ is weakly regular if and only if $Q_{A(\Lambda)}$ is acyclic.
        \item $\Lambda$ is Gorenstein if and only if $A(\Lambda)$ is a product of some finite dimensional non-simple self-injective radical square zero algebras if and only if each connected component of $Q_{A(\Lambda)}$ is a cycle with trivial valuations.
        \item $\Lambda$ is Iwanaga-Gorenstein if and only if $A(\Lambda)$ is Iwanaga-Gorenstein if and only if each connected component of $Q_{A(\Lambda)}$ is acyclic or is a cycle with trivial valuations.
        \item $\Lambda$ is sg-Hom-finite if and only if $A(\Lambda)$ is sg-Hom-finite if and only if $Q_{A(\Lambda)}$ is obtained from a disjoint union of oriented cycles with trivial valuations by adjoining sources and sinks with arbitrary values repeatedly.
    \end{enumerate}
\end{theorem}

This paper is organized as follows. In Section \ref{sec: Preliminaries}, we recall the theory of the stabilization of a left triangulated category and collect relevant facts about von Neumann regular rings and Gorenstein-projective modules. These will be used to study the singularity categories and to classify the specific classes of B{\"a}ckstr{\"o}m orders. Sections \ref{subsec: General Theory of Orders} and \ref{subsec: The Representation Theory of a B-Order} review the general theory of orders and the representation theory of B{\"a}ckstr{\"o}m orders developed by Ringel and Roggenkamp. Section \ref{sec:sing} is devoted to proving our main results, followed by the classification criteria. Finally in Section \ref{sec: example}, we provide some examples.

\begin{conventions and notations}\label{conventions}
    Throughout this article, all modules are considered as finitely generated right modules unless otherwise specified. For a noetherian ring $\Lambda$, we denote by $\mod \Lambda$ the category of finitely generated $\Lambda$-modules. Let $\mathcal{C}$ be an additive category. We denote by $\ind \mathcal{C}$ the isomorphism classes of the objects in $\mathcal{C}$. Let $M\in\mathcal{C}$. We denote by $\add M=\add_\mathcal{C} M$ the smallest subcategory of $\mathcal{C}$ containing $M$ which is closed under finite direct sums, direct summands and isomorphism. If $\mathcal{C}=\mod \Lambda$, we abbreviate $\add_{\mod \Lambda} M$ as $\add_\Lambda M$ and denote by $\proj \Lambda := \add_\Lambda \Lambda$ the category of finitely generated projective $\Lambda$-modules. Denote $\smod \Lambda$ as the stable module category of $\Lambda$ and $\Omega: \smod \Lambda \to \smod \Lambda$ as the syzygy functor. For $M\in \mod \Lambda$, $\sadd_\Lambda M$ is defined as $\add_{\smod\Lambda} M$. We drop the subscripts if the base ring is clear.
\end{conventions and notations}

\section{Preliminaries}\label{sec: Preliminaries}

\subsection{Stabilizations and Singularity Categories} \label{subsec: Stabilizations and Singularity Categories}

Let $\Lambda$ be a noetherian ring, $\rm{D}^{\rm{b}}(\Lambda)$ the bounded derived category of $\mod \Lambda$ and $\rm{K}^{\rm{b}}(\proj \Lambda)$ the bounded homotopy category of $\proj \Lambda$. The \textit{singularity category} $\Dsg(\Lambda)$ of $\Lambda$ is the Verdier quotient $\Dsg(\Lambda):=\rm{D}^{\rm{b}}(\Lambda)/\rm{K}^{\rm{b}}(\proj \Lambda)$. Denote by [1] the suspension functor of $\Dsg(\Lambda)$ and by $q: \rm{D}^{\rm{b}}(\Lambda) \to \Dsg(\Lambda)$ the quotient functor. The composition $\mod A \hookrightarrow \mathbf{\rm{D}}^{\rm{b}}(\Lambda) \xrightarrow{q} \Dsg(\Lambda)$ induces a functor $\smod A \to \Dsg(\Lambda)$ since projectives vanish, which is also denoted by $q$.

The so-called \textit{stabilization} is a useful tool to study singularity categories. It was first studied by Heller \cite{Hel68} in the context of topology and later studied in general settings by Keller and Vosseick \cite{KV87} and Beligiannis \cite{Bel00}. Before introducing stabilizations, we recall the notion of \textit{left triangulated categories}, which are one-sided versions of triangulated categories.  We refer to \cite[Section 2]{Che18} for an introduction to left triangulated categories and their stabilizations.

A \textit{left triangulated category} $(\mathcal{C},\Omega,\Delta)$ consists of an additive category $\mathcal{C}$, an additive endofunctor $\Omega$  called a \textit{loop functor} on $\mathcal{C}$, and a class $\Delta$ of sequences in $\mathcal{C}$ called \textit{left triangles} having the form $\Omega Z \to X \to Y \to Z$ which satisfy the left-sided version of the axioms of a triangulated category. In particular, every triangulated category is left triangulated. The notions of a \textit{left triangle functor} and a \textit{left triangulated subcategory} are defined similarly to those in the theory of triangulated categories; see \cite[Definition 2.4]{BM94}.

\begin{remark}
    It is clear that a left triangulated category is triangulated if and only if $\Omega$ is an autoequivalence. We call a left triangulated subcategory $(\mathcal{B}, \Omega)$ of $(\mathcal{C}, \Omega)$ a \textit{triangulated subcategory} of $(\mathcal{C}, \Omega)$ if $(\mathcal{B}, \Omega)$ is itself a triangulated category or equivalently $\Omega|_\mathcal{B}$ gives an autoequivalence of $\mathcal{B}$.
\end{remark}

Recall that an abelian category is called \textit{semisimple} if each short exact sequence splits. The (left) triangulated structure of a (left) triangulated category is called \textit{trivial} if each (left) triangle is isomorphic to a direct sum of some trivial (left) triangles. The following properties are elementary.

\begin{lemma}\cite[Lemma 3.4]{Che11}\label{lem:semisimpletriangle}
        Let $\mathcal{C}$ be a semisimple abelian category and $\Omega$ be an endofunctor (resp. autoequivalence) on $\mathcal{C}$. There exists exactly one left triangulated structure (resp. triangulated structure) on $\mathcal{C}$, i.e. the trivial structure, with respect to $\Omega$ as the loop functor (resp. cosuspension functor).
\end{lemma}

The \textit{stabilization} $(\mathcal{S}(\mathcal{C}),S)$ of a left triangulated category $(\mathcal{C},\Omega, \Delta)$ is a triangulated category which formally inverts the loop functor $\Omega$. It contains a triangulated category $\mathcal{S}(\mathcal{C})=(\mathcal{S}(\mathcal{C}), \tilde{\Omega}, \tilde{\Delta})$ where $\tilde{\Omega}$ is the cosuspension functor and a left triangle functor $S: \mathcal{C} \to \mathcal{S}(\mathcal{C})$ called \textit{the stabilization functor}. The construction of $(\mathcal{S}(\mathcal{C}),S)$ is as follows. The objects are $[M,m]$, where $M \in \mathcal{C}$ and $m\in \Z$. The morphism set of $[M,m],[N,n]\in\mathcal{S}(\mathcal{C})$ is 
\begin{equation*}
   \Hom_{\mathcal{S}(\mathcal{C})}([M,m],[N,n])=\underset{i\geq{m,n}}{\varinjlim} \Hom_\mathcal{C}(\Omega^{i-m}M,\Omega^{i-n}N).
\end{equation*}
The cosuspension functor $\tilde{\Omega}$ sends $[M,m]$ to $[M,m-1]$ acting in an obvious way on morphisms. Moreover, there is a natural isomorphism $[M,m-1]=\tilde{\Omega}[M,m]\simeq[\Omega(M),m]$.
The stabilization functor $S: \mathcal{C} \to \mathcal{S}(\mathcal{C})$ sends $M$ to $[M,0]$ on objects and $f:M \to N$ to the zero-representative  on morphisms. We remark that the stabilization functor is neither full nor faithful in general. 

The stabilization satisfies the following universal property: for any triangulated category $\mathcal{D}$ and any left triangle functor $F:\mathcal{C} \to \mathcal{D}$, there exists a unique triangle functor $\tilde{F}: \mathcal{S}(\mathcal{C}) \to \mathcal {D}$ such that $\tilde{F}S=F$. Consider a left triangle functor $F: (\mathcal{C},\Omega) \to (\mathcal{C'}, \Omega')$ between two left triangulated categories. By the universal property, it induces a triangle functor $\tilde{F}$ such that the following diagram commutes. 
\begin{equation*}
    \begin{CD}
        (\mathcal{C},\Omega) @>S>> S(\mathcal{C},\Omega) \\
        @VFVV @V\tilde{F}VV \\
        (\mathcal{C'},\Omega') @>S'>> S(\mathcal{C'},\Omega') 
    \end{CD}
\end{equation*}
The following lemma provides a criterion for when $\tilde{F}$ is an equivalence.

\begin{lemma}\cite[Proposition 2.2]{Che18}\label{pro:criterion for tilde{F}}
    The following statements hold.
    \begin{enumerate}[label=\textnormal{(\arabic*)}]
        \item $\tilde{F}$ is faithful if for any morphism $g: FX \to FY$, there 
        exists a morphism $f: \Omega^i(X) \to \Omega^i(Y)$ for some $i \geq 0$ such that $\Omega'^i(g)=\delta_Y^i \circ F(f) \circ (\delta_X^i)^{-1}$. Here $\delta: \Omega' \circ F \xrightarrow{\simeq} F \circ \Omega$ denotes the natural isomorphism of the left triangle functor $F$.
        \item $\tilde{F}$ is full if and only if for any two morphisms $f,f': X \to Y$ in $\mathcal{C}$ such that $F(f)=F(f')$, there exists $i\geq 0$ such that $\Omega^i(f)=\Omega^{i}(f')$.
        \item $\tilde{F}$ is dense if and only if for any $C' \in \mathcal{C'}$, there exist $i\geq 0$ and $X\in \mathcal{C}$ such that $\Omega'^{i}(C')\simeq F(X)$.
    \end{enumerate}
        In particular, if $(\mathcal{C}, \Omega)$ is a left triangulated subcategory of $(\mathcal{C'}, \Omega)$, $F$ is the inclusion functor and there exists $i\geq 0$ such that $\Omega^{i}(\mathcal{C'}) \subseteq \mathcal{C}$, then $\tilde{F}: (\mathcal{C}, \Omega) \xrightarrow{\simeq} (\mathcal{C'}, \Omega)$.
\end{lemma}

A typical example of a left triangulated category is the stable module category $\smod \Lambda$ of a noetherian ring $\Lambda$. The loop functor is given by the syzygy functor $\Omega: \smod \Lambda \to \smod \Lambda$. An exact sequence $0 \to L \to M \to N \to 0$ in $\mod \Lambda$ gives a left triangle $\Omega N \to L \to M \to N$ and any left triangle is isomorphic to a left triangle obtained in this way; see \cite[Proposition 2.10]{BM94}. For the stable module category $\smod \Lambda$, we denote its stabilization by $\mathcal{S}(\Lambda):=\mathcal{S}(\smod \Lambda, \Omega)$. The canonical functor $q: \smod \Lambda \to \Dsg(\Lambda)$ is in fact a left triangle functor. By the universal property of stabilizations, we have a triangle functor $\tilde{q}: (\mathcal{S}(\Lambda),\tilde{\Omega}^{-1}) \to (\Dsg(\Lambda),[1])$.

\begin{theorem}\cite[Theorem 3.8]{Bel00,KV87} \label{thm:stabilization}
    The functor $\tilde{q}: (\mathcal{S}(\Lambda),\tilde{\Omega}^{-1}) \to (\Dsg(\Lambda),[1])$ is a triangle equivalence. As a consequence, we have
    \begin{equation*}
        \Hom_{\Dsg(\Lambda)}(M,N)\simeq \underset{i\geq0}{\varinjlim} \sHom_\Lambda(\Omega^{i}M,\Omega^{i}N)
    \end{equation*}
    for any $M,N\in \mod\Lambda$.
\end{theorem}

The following proposition will be useful in the proof of our main theorem. It is included in \cite[Section 2]{Che11}. For the convenience of readers, we provide a brief proof here. Recall that an additive category is called \textit{idempotent complete} if each idempotent $e: X\to X$ splits, i.e. it admits a factorization $X\xrightarrow{u} Y \xrightarrow{v} X$ with $uv=\id_Y$. For example, a Krull-Schmidt category is idempotent complete. 

\begin{proposition}\label{prop:Dsg_iso_projEnd}
    Let $\Lambda$ be a noetherian ring. Assume that $\mod\Lambda$ is a Krull-Schmidt category and that there is a $\Lambda$-module $X$ and $n\geq0$ such that $\Omega^n(\smod\Lambda)\subseteq \sadd X$, where $\sadd X$ is the full subcategory of $\smod\Lambda$ generated by $X$. The following statements hold.
    \begin{enumerate}[label=\textnormal{(\arabic*)}]
        \item The singularity category $\Dsg(\Lambda)$ has an additive generator $q(X)$, i.e. $\Dsg(\Lambda)= \add q(X)$.
        \item The singularity category $\Dsg(\Lambda)$ is idempotent complete and we have
        \begin{equation*}
            \Hom_{\Dsg(\Lambda)}(q(X),?): \Dsg(\Lambda) \xrightarrow{\sim} \proj \End_{\Dsg(\Lambda)}(q(X)).
        \end{equation*}
    \end{enumerate}
    \begin{proof}
        Since $\mod \Lambda$ is Krull-Schmidt, $\smod \Lambda$ is also Krull-Schmidt and in particular is idempotent complete. It is easy to see the stabilization $\mathcal{S}(\Lambda)$ is idempotent complete by construction and so is $\Dsg(\Lambda) \simeq \mathcal{S}(\Lambda)$. Then the last statement of (2) follows from the projectivization of an idempotent complete category with an additive generator \cite[Proposition II.2.1]{ARS97}. 
        
        For (1), since $\Omega^n(\smod \Lambda)\subseteq \sadd E$, we have a sequence 
        \begin{equation*}
            \sadd X \supseteq \sadd (\Omega^n(\smod \Lambda))\supseteq \sadd (\Omega^{n+1}(\smod \Lambda)) \supseteq \cdots.
        \end{equation*}
        This sequence must terminate, because $\sadd X$ has only finitely many indecomposable modules up to isomorphism. Thus, we may assume without loss of generality that
        \begin{equation*}
            \sadd X \supseteq \sadd (\Omega^n(\smod \Lambda))=\sadd (\Omega^{n+1}(\smod \Lambda)).
        \end{equation*}

        Let $M[m]\in\Dsg(\Lambda)$. Since in the singularity category, the negative shift corresponds to the syzygy, we have $M[m]\simeq\Omega^{n_1}(M)[m+n_1]$ for any $n_1\geq-m,n$. Since $\sadd(\Omega^{n_1}(\smod \Lambda))=\sadd(\Omega^{n+m+n_1}(\smod\Lambda))$, then $M[m]\simeq\Omega^{n_1}(M)[m+n_1]$ is a direct summand of $\Omega^{n+m+n_1}(Y)[m+n_1]\simeq\Omega^n(Y)$ for some $Y \in \smod \Lambda$. By assumption, $\Omega^n(Y)\in \sadd X$. Then $M[m] \in \add q(X)$ as desired.
    \end{proof}
\end{proposition}

\subsection{Von Neumann Regular Rings}\label{sec:v-N alg}

The singularity category of a B{\"a}ckstr{\"o}m order will be described as the category of finitely generated projective modules over a certain von Neumann regular algebra. In this subsection, we review some fundamental properties of a von Neumann regular ring. We provide brief proofs for the reader's convenience.

\begin{thmdef}\cite[Theorem-Definition 11.24]{Fai73} \label{thmdef:von Neumann regular}
    Let $V$ be any ring. The following statements are equivalent.
    \begin{enumerate}[label=\textnormal{(\arabic*)}]
        \item Each right $V$-module is flat.
        \item For any $a\in V$, there exists $x\in V$ such that $a=axa$.
        \item Each finitely generated right ideal is generated by an idempotent.
    \end{enumerate}
    \rm{A ring $V$ satisfying one of the conditions above is called a \textit{von Neumann regular} ring.}
\end{thmdef}

\begin{proposition}\label{cor:semisimple is von Neumann regular} \label{cor:direct limit of von Neumann regular}\label{cor:proj of v-N is semisimple abelian}
    Let $V$ be any ring. Then the following statements hold.
    \begin{enumerate}[label=\textnormal{(\arabic*)}]
        \item If $V$ is an artinian semisimple ring, then $V$ is von Neumann regular.
        \item A direct limit of von Neumann regular rings is also von Neumann regular. In particular, a direct limit of artinian semisimple rings is von Neumann regular. 
        \item If $V$ is a von Neumann regular ring, then the category $\proj V$ is semisimple abelian.
    \end{enumerate}
    \begin{proof}
        (1) and (2) follow immediately from Theorem-Definition \ref{thmdef:von Neumann regular} (1) and (2). For (3), we only need to show $\proj V$ is abelian. To show this, it is enough to prove that $\proj V$ is closed under kernels and cokernels. Let $f:P \to Q$ be a homomorphism for some $P,Q\in\proj V$. Then we have an exact sequence $P \xrightarrow{f} Q \to \Coker f\to 0$,
        which implies $\Coker f$ is finitely presented. By the condition Theorem-Definition \ref{thmdef:von Neumann regular} (1), we have $\Coker f\in \proj V$. Here we use the well-known fact that a module over any ring is finitely generated and projective if and only if it is finitely presented and flat. So, $\Im f \in \proj V$, which implies $\Ker f$ is also in $\proj V$.
    \end{proof}
\end{proposition}

\begin{defprop} \label{ex:V(D,M)}
    Assume that $D$ is a semisimple ring and $M$ is a $(D,D)$-bimodule. Let $M^{\otimes 0}:=D$ and $M^{\otimes i}:=M \otimes_D M^{\otimes{i-1}}$ for $i \geq 1$. 
    Define
    \begin{equation*}
        V=V(D,M):=\underset{i\geq0}{\varinjlim} \End_D(M^{\otimes i}) \quad \text{and} \quad K=K(D,M):=\underset{i\geq1}{\varinjlim} \Hom_D(M^{\otimes i}, M^{\otimes i-1}).
    \end{equation*}
    Then 
    \begin{enumerate}[label=\textnormal{(\arabic*)}]
        \item $V$ is a von Neumann regular algebra and $K$ is a $(V,V)$-bimodule.
        \item The category $\proj V$ is semisimple abelian and $(\proj V, ? \otimes_V K)$ is a trivial left triangulated category.
    \end{enumerate}
    \begin{proof}
        These follow by definition, Proposition \ref{cor:semisimple is von Neumann regular} and Lemma \ref{lem:semisimpletriangle}.
    \end{proof}
\end{defprop}

We review Xiao-wu Chen's description of the singularity category of a finite dimensional radical square zero algebra. 

\begin{theorem}\cite[Theorem 3.8]{Che11}\label{thm:DsgArtin}
    Let $A$ be a finite dimensional radical square zero algebra and let $J$ be the radical of $A$. Then $\Omega(\smod A)\subseteq \sadd (A/J)$ and the syzygy functor $\Omega$ is isomorphic to 
    \begin{equation*}
        \Omega \simeq ?\otimes_{A/J} J: \sadd (A/J) \to \sadd (A/J).
    \end{equation*}
    Define 
    \begin{equation*}
        V(A):=V(A/J,J)=\underset{i\geq0}{\varinjlim} \End_{A/J}(J^{\otimes i}) \quad {\rm{and}} \quad K(A):=K(A/J,J)=\underset{i\geq 1}{\varinjlim} \Hom_{A/J}(J^{\otimes i}, J^{\otimes i-1}).
    \end{equation*}
    Then $V(A)$ is a von Neumann regular algebra and there is a triangle equivalence 
    \begin{equation*}
        (\Dsg(A),[1]) \simeq (\proj V(A), ? \otimes_{V(A)}K(A)),
    \end{equation*}
    where the triangulated structures are trivial.
\end{theorem}

\subsection{Gorenstein-projective Modules}\label{subsec:Gorenstein-projective modules}

In this subsection, we collect relevant facts about Iwanaga-Gorenstein rings and Gorenstein-projective modules, which will be used to classify Iwanaga-Gorenstein and Gorenstein B{\"a}ckstr{\"o}m orders in Section \ref{subsec: Classifications of some Classes of B-Orders}. Throughout this subsection, we assume that $\Lambda$ is a noetherian ring.

\begin{definition}
    A noetherian ring $\Lambda$ is called \textit{Iwanaga-Gorenstein} if $\id_\Lambda \Lambda< \infty$ and $\id_{\Lambda^{\op}} \Lambda<\infty$. In this case, $\id_\Lambda \Lambda=\id_{\Lambda^{\op}} \Lambda<\infty$; see \cite[Proposition 9.1.8]{EJ00}.
\end{definition}

\begin{definition}
    A finitely generated $\Lambda$-module $M$ is called \textit{Gorenstein-projective} if there is an exact sequence
        $P_\bullet=(\cdots \to P_{1} \xrightarrow{f_{_1}} P_{0} \xrightarrow{f_{0}} P_{-1} \xrightarrow{f_{-1}} P_{-2} \to \cdots)$,
    where $P_i\in\proj\Lambda$ and $M=\Ker f_{-1}$, such that $\Hom_\Lambda(P_\bullet,Q)$ is also an exact sequence for any $Q\in\proj\Lambda$.

    We denote the full subcategory consisting of Gorenstein-projective modules by $\Gp \Lambda$.
\end{definition}

It is clear by definition that $\proj\Lambda \subseteq \Gp \Lambda$, $\Omega(\Gp \Lambda)\subseteq\Gp \Lambda$, $\Gp \Lambda\subseteq \bigcap_{i=1}^{\infty}\Omega^i(\mod \Lambda)$ and $\pd M=\infty$ for any $M\in \Gp \Lambda \backslash \proj \Lambda$. Moreover, $\Gp \Lambda$ is a Frobenius category with projective-injective objects $\proj \Lambda$ and thus the stable category $\sGp \Lambda$ is a triangulated subcategory of $s\mod \Lambda$ \cite{Hap88}. Moreover, Beligiannis proved that $\sGp\Lambda$ is the largest triangulated subcategory of $\smod \Lambda$; see the following proposition.

\begin{proposition}\cite[Proposition 2.13]{Bel00}\label{prop: maximum property of GP}
    If a full additive subcategory $\mathcal{C}$ of $\smod \Lambda$ is a triangulated subcategory of $\smod \Lambda$, then $\mathcal{C}\subseteq\sGp\Lambda$.
\end{proposition}

As for projective modules, we can define Gorenstein-projective dimension for Gorenstein-projective modules.

\begin{definition}
    A finitely generated module $M$ is said to have \textit{Gorenstein-projective dimension} $\Gpd M$ at most $n$ if there exists a finite exact sequence $0 \to G_n \to \dots \to G_0 \to M \to 0$ with $G_i\in\Gp \Lambda$. If no such a finite exact sequence exists, then $\Gpd M$ is defined by $\infty$.
\end{definition}

Consider the canonical triangle functor $i: \sGp\Lambda \hookrightarrow \smod \Lambda \xrightarrow{S} \mathcal{S}(\Lambda)\simeq\Dsg(\Lambda)$. The property that each module has finite Gorenstein-projective dimension can be characterized by the requirement that this functor is an equivalence.

\begin{theorem}\label{thm: criterion for IG}
    Let $\Lambda$ be a noetherian ring. The canonical functor $i:\Gp\Lambda \to \Dsg(\Lambda)$ is full and faithful. Moreover, this functor is dense if and only if $\Gpd M <\infty$ for any $M\in\mod\Lambda$.
    \begin{proof}
        The sufficiency is by Buchweitz and Happel \cite{Buc97,Hap91}. The necessity is by Beligiannis \cite[Theorem 3.8]{Bel00}.
    \end{proof}
\end{theorem}

Any module over an Iwanaga-Gorenstein ring has finite Gorenstein-projective dimension; see the following proposition. The converse is not true in general, as there need not exist a uniform bound on Gorenstein-projective dimension. However, the converse does hold for an order; see Theorem \ref{thm:criterion for Gor order}.

\begin{proposition} \label{prop: IG implies finite Gpd} \cite[Corollary 11.5.3]{EJ00}
    Let $\Lambda$ be an Iwanaga-Gorenstein ring. Then $\Gpd M <\infty$ for any $M\in\mod\Lambda$. In particular, the canonical functor $i:\Gp\Lambda \to \Dsg(\Lambda)$ is an equivalence. 
\end{proposition}

We end this subsection with a further property of Gorenstein-projective dimension. We remark that the notions \textit{Gorenstein-projective objects} and \textit{Gorenstein-projective dimension} can be defined similarly for an exact category with enough projectives.

\begin{lemma}\label{lem:lem to Gor}
    Let $\mathcal{C}$ be an exact category with enough projectives and injectives. Then if $M\in \mathcal{C}$ satisfies $\id_\mathcal{C}M<\infty$, then $\Gpd_\mathcal{C} M=\pd_{\mathcal{C}} M$. 
    \begin{proof}
        This property was proved by Holm in the setting of module categories \cite[Theorem 2.2]{Hol04b}. The argument is still valid for arbitary exact categories. 
    \end{proof}
\end{lemma}

\section{The Representation Theory of B{\"a}ckstr{\"o}m Orders}

\subsection{General Theory of Orders}\label{subsec: General Theory of Orders}
In this subsection, we review the general theory for orders. Our main references are \cite[Chapter 3]{CR90}, \cite[Section 1]{HN94} and \cite[Chapter 39]{Rei03}. For the reader's convenience, we include short proofs.

Throughout this subsection we assume that $R$ is a complete regular local ring of Krull dimension $d$ and $K$ is the fraction field of $R$.

\begin{definition}\label{def:order}
    An \textit{$R$-order} $\Lambda$ is an $R$-algebra satisfying $\Lambda$ is free of finite rank as an $R$-module.
\end{definition}

For an order of positive Krull dimension, we mainly study the category of Cohen-Macaulay modules instead of the whole module category.

\begin{definition}
    A finitely generated $\Lambda$-module $M$ is called a \textit{Cohen-Macaulay $\Lambda$-module} if $M$ is free as an $R$-module. The full subcategory of $\mod \Lambda$ consisting of Cohen-Macaulay $\Lambda$-modules is denoted by $\CM\Lambda$.
\end{definition}

\begin{remark}
    The category $\CM\Lambda$ is an extension-closed subcategory of $\mod\Lambda$ and hence an exact category. We will see that it has enough projectives and injectives.
\end{remark}

\begin{remark}\label{rmk: Krull-Schmidt Property}
    Since $R$ is a complete local ring, the category $\mod \Lambda$ is Krull-Schmidt and so is $\CM \Lambda$; see \cite[Section 1.2]{LW12}.
\end{remark}

\begin{remark}
    In some references, a Cohen-Macaulay $\Lambda$-module $M$ is called a \textit{$\Lambda$-lattice} since it is free over $R$. The terminology of Cohen-Macaulay $\Lambda$-module arises from the fact that the defining condition is equivalent to $M$ being a (maximal) Cohen-Macaulay $R$-module (i.e. $\depth_R M=\dim R$); see \cite[Theorem 1.3.3]{BH93}. Moreover, when $\dim R=1$, the condition is further equivalent to $M$ being $R$-torsion-free; see \cite[Section 12.1]{DF04}.
\end{remark}

\begin{definition}
    Let $\Lambda$ be an $R$-order. 
    \begin{enumerate}
        \item Assume that $M,N\in\CM\Lambda$. Then $M$ is called an \textit{overmodule} of $N$ if $N\subseteq M \subseteq N \otimes_R K$ and $N$ is a $\Lambda$-submodule of $M$.
        \item An \textit{overorder} $\Gamma$ of $\Lambda$ is an $R$-order such that $\Lambda$ is an $R$-subalgebra of $\Gamma$ and $\Gamma$ is a $\Lambda$-overmodule of $\Lambda$. 
        \item $\Lambda$ is called a \textit{maximal order} if $\Lambda$ has no proper overorder.
        \item $\Lambda$ is called a \textit{regular order} if $\gd\Lambda=\dim R$. If $\dim R=1$, a regular order is also called a \textit{hereditary order}.
    \end{enumerate}
\end{definition}

The following properties are elementary for orders.

\begin{proposition}\label{prop:hereditary order} \label{prop:duality}
    Let $\Lambda$ be an $R$-order. 
    \begin{enumerate}[label=\textnormal{(\arabic*)}]
        \item There is an exact duality $(?)^{\vee}:=\Hom_R(?,R): \CM \Lambda \xrightarrow{\simeq} \CM \Lambda^{\op}$. 
        \item Define $\inj \Lambda:=\add (_\Lambda\Lambda)^{\vee}\subseteq\CM\Lambda$. Then the duality in (1) gives dualities $(?)^{\vee}: \proj \Lambda \xrightarrow{\simeq} \inj \Lambda^{\op}$ and $(?)^{\vee}: \inj \Lambda \xrightarrow{\simeq} \proj \Lambda^{\op}$.
        \item $\proj \Lambda=\{\text{projective objects in } \CM\Lambda\}$ and $\inj \Lambda=\{\text{injective objects in } \CM\Lambda\}$. Moreover, $\CM\Lambda$ has enough projectives and injectives.
        \item  $\Lambda$ is regular if and only if $\CM \Lambda = \proj \Lambda$. 
    \end{enumerate}
     \begin{proof}
        (1) follows by definition and $(?):\proj R \xrightarrow{\simeq} \proj R$ being an exact duality. (2)(3) follow easily by (1). For (4), see \cite[Proposition 2.17]{IW14}.
    \end{proof}
\end{proposition}

\begin{lemma}\label{lem: compute Hom for orders}
    Let $\Lambda$ be an R-order and $M,N\in \CM \Lambda$. Then 
    \begin{enumerate}[label=\textnormal{(\arabic*)}]
        \item $(M)^{\vee}\simeq\{f\in \Hom_{K}(M\otimes_R K,K) \mid f(M)\subseteq R\}$
        \item $\Hom_\Lambda(M,N)\simeq\{f\in \Hom_{\Lambda\otimes_R K}(M\otimes_R K,N\otimes_R K) \mid f(M)\subseteq N\}$.
        \item Let $M$ be an overmodule of $N$. Then $N\otimes_R K=M\otimes_R K$. Moreover $(N)^{\vee}$ is an overmodule of $(M)^{\vee}$ in $\CM \Lambda^{\op}$.
    \end{enumerate}
    \begin{proof}
        (2): The left-hand side is an $R$-submodule of $\Hom_{\Lambda\otimes_R K}(M\otimes_R K,N\otimes_R K)$ since $N$ is $R$-free. Then the equality can be checked easily. The proof of (1) is similar. (3) follows immediately by (1).
    \end{proof}
\end{lemma}

\begin{proposition}\label{prop: overmodule} \label{prop:overorder} 
    Let $\Gamma$ be an overorder of an $R$-order $\Lambda$. Then the restriction functor 
    \begin{equation*}
        \res: \CM\Gamma \to \CM \Lambda
    \end{equation*}
     is fully faithful.
    \begin{proof}
        This proposition follows easily from the definition and Lemma \ref{lem: compute Hom for orders}
    \end{proof}
\end{proposition}

For $R$-orders, being Iwanaga-Gorenstein is equivalent to the property that each module has finite Gorenstein-projective dimension (c.f. Proposition \ref{prop: IG implies finite Gpd}). This is a natural generalization of the artinian case \cite{AB69}.

\begin{lemma}\label{lem: R order Gpd finite implies IG}
     Let $\Lambda$ be an $R$-order. If $\Gpd M < \infty$ for any $M\in\mod\Lambda$, then $\Lambda$ is Iwanaga-Gorenstein.
     \begin{proof}
        Recall that our modules are right modules and $\inj \Lambda$ is the subcategory of the injectives in $\CM\Lambda$.
        Firstly, we show $\id_{\Lambda^{\op}} \Lambda < \infty$. We have $\Gpd_\Lambda (\Lambda^{\vee}) < \infty$.  Since $\Lambda^{\vee} \in \inj \Lambda$, by Lemma \ref{lem:lem to Gor}(2), we have $\pd_\Lambda (\Lambda^{\vee}) = \Gpd_\Lambda (\Lambda^{\vee}) < \infty$. So, $\id_{\CM\Lambda^{\op}} \Lambda < \infty$. By \cite[1.2(1)]{Iya05a}, each injective object in $\CM\Lambda$ has injective dimension $d$ in $\mod\Lambda$. Thus, $\id_{\Lambda^{\op}} \Lambda =\id_{\mod \Lambda^{\op}} \Lambda < \infty$.

        Secondly, we show $\id_\Lambda \Lambda < \infty$. Let $S:=\Lambda/\rad\Lambda$. Since $S':=\Omega^d(S) \in \CM\Lambda$, we have $\Gpd_\Lambda S'=n<\infty$ for some $n\geq 1$. Then by \cite[Theorem 2.20]{Hol04a}, $\Omega^n(S')\in \Gp \Lambda$. Thus, $\Ext_\Lambda^i(S,\Lambda)\simeq \Ext_\Lambda^{i-n-d}(\Omega^{n+d}(S),\Lambda)=\Ext_\Lambda^{i-n-d}(\Omega^{n}(S'),\Lambda)=0$ for $i\geq n+d+1$, which implies $\id_\Lambda\Lambda\leq n+d <\infty$ by the Bass-type lemma; see Corollary \cite[Corollary 3.5(3)]{GN02}.
     \end{proof}
\end{lemma}

\begin{theorem}\label{thm:criterion for Gor order}
     Let $\Lambda$ be an $R$-order. The following statements are equivalent.
    \begin{enumerate}[label=\textnormal{(\arabic*)}]
        \item $\Lambda$ is Iwanaga-Gorenstein .
        \item $\Gpd M < \infty$, for any $M\in\CM\Lambda$.
        \item The canonical functor $\sGp \Lambda \to \Dsg(\Lambda)$ is an equivalence.
    \end{enumerate}
    \begin{proof}
        (1) $\Rightarrow$ (2) is by Proposition \ref{prop: IG implies finite Gpd}. 
        (2) $\Rightarrow$ (1) is by Lemma \ref{lem: R order Gpd finite implies IG}.
        (2) $\Leftrightarrow$ (3) is by Theorem \ref{thm: criterion for IG}.
    \end{proof}
\end{theorem}


In the rest of this subsection, we assume that $R$ is a complete discrete valuation ring.

\begin{defprop}\label{def: MGamma}
    Let $\Gamma$ be an overorder of $\Lambda$. For $M\in\CM\Lambda$, define 
    \begin{equation*}
        M\Gamma:=\{\sum_i (m_i\otimes 1) \gamma_i \in M\otimes_R K \mid m_i \in M, \gamma_i \in \Gamma\}\in \CM \Gamma.  
    \end{equation*}
    This induces a functor $(?)\Gamma: \CM\Lambda \to \CM\Gamma$. Moreover, there exist adjoint pairs $((?)\Gamma, \res)$ and $(\res, \Hom_\Lambda(\Gamma,?))$.
    \begin{proof}
        Since $M\Gamma$ is $R$-torsion-free and $\dim R=1$, $M\Gamma\in\CM\Gamma$. The adjoint pairs can be checked easily.
    \end{proof}
\end{defprop}

There is a dual concept called \textit{coradical} to the radical of a Cohen-Macaulay module. 

\begin{defprop}\label{def:corad}
    For $M\in\CM\Lambda$, the \textit{coradical} $\corad M$ of $M$ is defined as 
    \begin{equation*}
        M\subseteq\corad M:=(\rad_{\Lambda^{\op}} (M^{\vee}))^\vee = \sum \{M\subsetneq N\mid \text{$N$ is a minimal overmodule of M}\}\in\CM\Lambda.
    \end{equation*}
    In particular, $\corad I$ is the unique minimal overmodule of $I$ if $I\in\inj\Lambda$ is indecomposable.
    \begin{proof}
        Note that $\CM\Lambda$ is closed under submodules since $\dim R=1$. Then the statements follow from the fact that $\rad M$ is the intersection of all maximal submodule of $M$ and from Propositions \ref{prop:duality} and \ref{prop: overmodule}.
    \end{proof}
\end{defprop}

The following is the structure theorem for hereditary orders which is analogous to Artin-Wedderburn theorem for finite dimensional semisimple algebras.

\begin{theorem}\label{thm:structure theorem of hereditary orders} \cite[Theorem 39.14]{Rei03} \cite[Theorem 1.7.1]{HN94}
    The following statements hold.
    \begin{enumerate}[label=\textnormal{(\arabic*)}]
        \item Let $D$ be a finite dimensional division $K$-algebra. Then there exists a unique maximal $R$-order $\Delta$ in $D$. Moreover, $\Delta$ is local and a (non-commutative) PID.
        \item Let $D$ and $\Delta$ be as in {\rm{(1)}}. Then $M_n(\Delta)\subseteq M_n(D)$ is a maximal order in $M_n(D)$. If $\Lambda$ is another maximal order, then $\Lambda$ is conjugate to $M_n(D)$ by an invertible element in $M_n(D)$.
        \item Assume that $\Lambda$ is a ring-indecomposable hereditary order, then $B:=\Lambda \otimes_R K$ is a simple $K$-algebra. Moreover, there exists an identification $B=M_n(D)$ for some finite dimensional division $K$-algebra $D$ such that
        \begin{equation*}
            \Lambda=
            \begin{pmatrix}
                \begin{matrix} M_{n_1}(\Delta) & \\ & M_{n_2}(\Delta) \end{matrix} & & \Delta \\
                    &   \ddots  & \\
                \rad \Delta & & M_{n_r}(\Delta)
            \end{pmatrix}_{n\times n}
        \end{equation*}
        where $n_1,n_2,\dots,n_r\geq 1$ and $\Delta$ is the unique maximal order in $D$.
    \end{enumerate}
\end{theorem}

\subsection{B{\"a}ckstr{\"o}m Orders and the Classical Result}\label{subsec: The Representation Theory of a B-Order}
In this subsection, we introduce the definition of B{\"a}ckstr{\"o}m orders and their representation theory developed by Ringel and Roggenkamp.
\begin{assumption}
    In the rest of this article, we assume that $(R,\mathfrak{m}=\pi R, k)$ is a complete discrete valuation ring and that $K=\Frac R$ is the fraction field of $R$.
\end{assumption} 

\begin{definition}\label{def:Backorder}
    A \textit{B{\"a}ckstr{\"o}m order} $\Lambda=(\Lambda,\Gamma)$ is an $R$-order $\Lambda$ together with an overorder $\Gamma$ of $\Lambda$ satisfying both of the following conditions
    \begin{enumerate}
        \item $\Gamma$ is a hereditary order;
        \item $\rad \Lambda =\rad \Gamma$. 
    \end{enumerate}
\end{definition}

\begin{remark}\label{rmk: CM Gamma=proj}
    By definition and Propositions \ref{prop:hereditary order} and \ref{prop:overorder}, we have $\CM \Gamma = \proj \Gamma \subseteq \CM \Lambda$.
\end{remark}

\begin{remark}
    Since $\Gamma$ is hereditary, $\Lambda \otimes_R K$ is semisimple by Theorem \ref{thm:structure theorem of hereditary orders} (3). It is equivalent to say that $\Lambda$ is an \textit{isolated singularity}, which is also equivalent to $\CM \Lambda$ admitting Auslander-Reiten sequences; see \cite{Aus86}.
\end{remark}

\begin{lemma}
    Let $\Lambda=(\Lambda,\Gamma)$ be a B{\"a}ckstr{\"o}m order. The following statements hold.
    \begin{enumerate}[label=\textnormal{(\arabic*)}]
        \item  The class of B{\"a}ckstr{\"o}m orders is invariant under Morita equivalence. Therefore, we may always assume that B{\"a}ckstr{\"o}m orders are basic.
        \item The hereditary order $\Gamma$ is uniquely determined by $\Lambda$. Therefore, for a B{\"a}ckstr{\"o}m order $\Lambda=(\Lambda,\Gamma)$, we simply write $\Lambda$ when no ambiguity arises. 
        \item For any $M\in\CM\Lambda$, the radical of $M$ coincides with the radical of $M\Gamma$ as a $\Lambda$-module, and also with the radical of $M\Gamma$ as a $\Gamma$-module. That is, we have 
    \begin{equation*}
        \rad_\Lambda M=\rad_\Lambda (M\Gamma)=\rad_\Gamma(M\Gamma)\in\CM\Gamma.
    \end{equation*}
    Therefore, we simply denote the radicals above by $\rad M$. \label{prop:rad GammaM=rad M}
    \end{enumerate}
    \begin{proof}
        For (1), see \cite[Lemma 2.2]{RR79}. For (2), there is a fact that, for any $R$-order $\Gamma$, $\Gamma$ is hereditary if and only if 
        \begin{equation*}
            \Gamma=\{a \in \Gamma \otimes_R K \mid (\rad\Gamma)a\subseteq\rad\Gamma\};
        \end{equation*}
        see \cite[Theorem 1.6]{HN94}. In our case, since $\rad\Lambda=\rad\Gamma$, we have 
        \begin{equation*}
            \Gamma=\{a\in \Lambda \otimes_R K \mid (\rad\Lambda)a\subseteq\rad\Lambda\},
        \end{equation*}
        which is uniquely determined by $\Lambda$. 
        For (3), since $\rad \Lambda =\rad \Gamma$ is a two-sided ideal both in $\Lambda$ and $\Gamma$, we have $\rad_\Gamma(M\Gamma)=(M\Gamma)(\rad\Gamma)=(M\Gamma)(\rad\Lambda)=M(\rad\Lambda)=\rad_\Lambda M$
        and the middle term is $\rad_\Lambda (M\Gamma)$.
    \end{proof}
\end{lemma}

From now on, we assume that $\Lambda=(\Lambda,\Gamma)$ is a B{\"a}ckstr{\"o}m order. Define a finite dimensional $k$-algebra $H$ associated to $\Lambda$ as follows 
\begin{equation}\label{eqa: H}
    H=H(\Lambda,\Gamma):=
    \begin{pmatrix}
        \Gamma / \rad \Gamma & 0 \\
        \Gamma / \rad \Gamma & \Lambda / \rad \Lambda 
    \end{pmatrix}
\end{equation}
This is a hereditary radical square zero algebra.

Define a functor 
\begin{equation*}
    \begin{tikzcd}[row sep=0.1em, column sep=3em]
        \mathbb{F}: \CM \Lambda \arrow[r] &  \mod H \\
            \quad M \arrow[r,|->,to path={([xshift=1em]\tikztostart.east) -- ([xshift=-0.25em]\tikztotarget.west)}] & ( M \Gamma /\rad M, M/\rad M)
    \end{tikzcd}
\end{equation*}
The $H$-action on $(M \Gamma /\rad M, M/\rad M)$ is given by the usual matrix multiplication and the canonical map $ M/\rad M \otimes_{\Lambda/\rad\Lambda} \Gamma/\rad \Gamma \to M\Gamma/\rad M$. For a homomorphism $f:M \to N$ in $\CM \Lambda$, it induces a $\Gamma$-module homomorphism $f_1=(f)\Gamma:  M\Gamma  \to  N\Gamma $ such that the following diagram commutes, in which the horizontal maps are the counits of the adjoint pair $((?)\Gamma, \res)$; see Definition-Proposition \ref{def: MGamma}
\begin{equation*}
    \begin{tikzcd}
        M \arrow[r,hook] \arrow[d,"f"] &  M\Gamma  \arrow[d,"f_1"]\\
        N \arrow[r,hook] & N\Gamma 
    \end{tikzcd}
\end{equation*}
Passing to the quotient, we have another commutative diagram:
\begin{equation*}
    \begin{tikzcd}
        M/\rad M  \arrow[r,hook] \arrow[d,"\bar{f}"] &  M\Gamma  /\rad M \arrow[d,"\bar{f_1}"]\\
        N/\rad N  \arrow[r,hook] &  N\Gamma /\rad N
    \end{tikzcd}
\end{equation*}
This defines an $H$-module homomorphism 
\begin{equation*}
    \mathbb{F}(f):=(\bar{f_1},\bar{f}):(M\Gamma  /\rad M, M/\rad M) \to (N\Gamma  /\rad N, N/\rad N). 
\end{equation*}

Denote by $\mod_s H$ the full subcategory of $\mod H$ whose objects contain no non-zero simple direct summand and $\simple H:=\ind (\mod H/\rad H)$. The next theorem by Ringel and Roggenkamp establishes the representation theory of B{\"a}ckstr{\"o}m orders.

\begin{theorem}\cite[Theorem 1.2]{RR79}\label{thm:rep_of_Back}
    Let $\Lambda=(\Lambda,\Gamma)$ be a B{\"a}ckstr{\"o}m order and keep the notations above. The functor 
    \begin{equation*}
        \begin{tikzcd}[row sep=0.1em, column sep=2em]
            \mathbb{F}: \CM \Lambda \arrow[r] &  \mod H \\
        \end{tikzcd}
    \end{equation*}
    is full with essential image $\mod_s H$. The kernel $\Ker \mathbb{F}_{M,N}$ for $M,N\in \CM\Lambda$ is $\Hom_\Lambda(M,\rad N)$. Moreover, this functor gives one-to-one correspondences:
    \begin{enumerate}[label=\textnormal{(\arabic*)}]
        \item $\ind(\CM \Lambda) \xleftrightarrow{1:1} \ind(\mod H)\backslash \simple H$.
        \item $\ind(\proj \Lambda) \xleftrightarrow{1:1} \ind(\proj H)\backslash\simple H$.
        \item $\ind(\proj \Gamma) \xleftrightarrow{1:1} \ind(\inj H)\backslash\simple H$.
    \end{enumerate}
\end{theorem}

\begin{remark}
    For the finite dimensional $k$-algebra $H$, we may associate a \textit{species} to it. The B{\"a}ckstr{\"o}m orders of \textit{finite Cohen-Macaulay type}, i.e. those with finitely many indecomposable Cohen-Macaulay modules up to isomorphism, are classified by the shape of the species \cite[Theorem III]{RR79}. More precisely, a B{\"a}ckstr{\"o}m order is of finite Cohen-Macaulay type if and only if its species is a finite union of Dynkin diagrams. Moreover, Roggenkamp also studied the Auslander-Reiten quiver of a B{\"a}ckstr{\"o}m order \cite{Rog83}. It is obtained by the Auslander-Reiten quiver of $\mod H$ under some combinatorial constructions.
\end{remark}

\subsection{Further Properties of B{\"a}ckstr{\"o}m orders}\label{subsec: The Representation Equivalence}
In this subsection, we investigate further properties of B{\"a}ckstr{\"o}m orders. In this subsection, we assume that $\Lambda=(\Lambda,\Gamma)$ is a B{\"a}ckstr{\"o}m order.

Let $N\in\CM\Lambda$. For any surjective homomorphism $f: P \rightarrowdbl N$ with $P\in \proj \Lambda$, it induces a commutative diagram of exact sequences:
\begin{equation}\label{cartesian_diagram}
            \begin{CD}
                        @.                    @.      0         @.          0           @.          \\  
                @.              @.                    @VVV                  @VVV                @.  \\        
                0       @>>>   \Omega(N)             @>>>    \rad P    @>>>        \rad N      @>>>    0   \\
                @.           @VV\rotatebox{270}{$\simeq$}V  @VVV             @VVV                @.  \\
                0       @>>>    \Omega(N)             @>>>    P         @>f>>       N           @>>>    0   \\
                @.              @.                    @VVV                  @VVV                @.  \\
                        @.                    @.      P/\rad P  @>\bar{f}>\simeq>   N/\rad N    @.       \\
                @.              @.                    @VVV                  @VVV                @.  \\
                        @.                    @.      0         @.          0           @.          
            \end{CD}
 \end{equation}

\begin{proposition}\label{prop:OmegaCMLambda}
    For a B{\"a}ckstr{\"o}m order $\Lambda$, we have $\Omega(\sCM\Lambda)\subseteq \sadd \Gamma$.
    \begin{proof}
        Let $N\in\CM\Lambda$. Since $\rad P, \rad N \in \CM \Gamma$ in the commutative diagram (\ref{cartesian_diagram}), then $\Omega(N) \in \CM \Gamma=\proj \Gamma$.
    \end{proof}
\end{proposition}

\begin{lemma}\label{lem:P(M,N)=(M,rad N)}
    The following statements hold.
    \begin{enumerate}[label=\textnormal{(\arabic*)}]
        \item For any $N \in \CM \Lambda$, the inclusion map $\rad N \hookrightarrow N$ factors through a projective $\Lambda$-module.
        \item Let $M,N \in \CM \Lambda$ be indecomposable and non-projective. Then we have
        \begin{equation*}
            \mathcal{P}(M,N)=\Hom_\Lambda(M,\rad N),
        \end{equation*} 
        where $\mathcal{P}(M,N)$ is the set of homomorphisms from $M$ to $N$ factoring through a projective module.
    \end{enumerate}
    \begin{proof}
        For (1), the exact sequence $0\to\Omega(N)\to\rad P \to \rad N \to 0$ splits since it is in $\CM\Gamma=\proj\Gamma$. So, $\rad N \hookrightarrow N$ factors through $P\in \proj \Lambda$.

        For (2), assume that $f: M \to N$ factors through $g: M \to P$ for some $P \in \proj \Lambda$. Since $M$ is not projective, then $\Im g \subseteq \rad P$. So, $\Im f \subseteq \rad N$. Conversely, assume that the image of $f: M \to N$ is in $\rad N$. This means $f$ factors through the inclusion map $\rad N \hookrightarrow N$. By (1), this inclusion factors through a projective $\Lambda$-module and so does $f$.
        \end{proof}
\end{lemma}

The following result shows that the functor $\mathbb{F}:\CM\Lambda \to \mod H$ in Theorem \ref{thm:rep_of_Back} induces a stable equivalence. This is a 1-dimensional analogue of a result for finite dimensional radical square zero algebras \cite[Theorem X.2.4]{ARS97}, providing evidence for our Slogan \ref{slogan}.

\begin{theorem}\label{thm:stable equivalence}
    The functor $\mathbb{F}: \CM \Lambda \to \mod H$ induces an equivalence
    \begin{equation*}
        \underline{\mathbb{F}}: \sCM \Lambda \xlongrightarrow{\sim} \smod_s H,
    \end{equation*}
    where $\smod_s H$ is the full subcategory of $\smod H$ whose objects do not contain any non-zero simple direct summand.
    \begin{proof}
        By construction, $\mathbb{F}$ sends $\proj \Lambda$ to $\proj H$. So it induces a functor $\underline{\mathbb{F}}: \sCM \Lambda \to \smod H$. It is full and dense in $\smod_s H$ by Theorem \ref{thm:rep_of_Back}. By Lemma \ref{lem:P(M,N)=(M,rad N)}, if $M,N$ have no projective direct summands, then $\Ker \mathbb{F}_{M,N}=\Hom_\Lambda(M,\rad N)=\mathcal{P}(M,N)$. Therefore, $\underline{\mathbb{F}}$ is also faithful.
    \end{proof}
\end{theorem}

Next, we analyze the connection between the syzygies of $N\in \proj \Lambda$ and of $\mathbb{F}(N) \in \mod_s H$ under the functor $\mathbb{F}: \CM \Lambda \to \mod H$.

\begin{proposition}\label{thm:syzygy_in_addGamma}
    Consider $N\in\CM\Lambda$ and $\mathbb{F}(N) \in \mod_s H$. The following statements hold.
    \begin{enumerate}[label=\textnormal{(\arabic*)}]
        \item Assume that $f:P \rightarrowdbl N$ is a projective cover of $N$ with $P\in\proj \Lambda$. It gives rise to a commutative diagram \rm{(\ref{CD:big_comm_diagram})} of exact sequences.
        \item Under the assumptions in (1), the homomorphism $\mathbb{F}(f):\mathbb{F}(P) \rightarrowdbl \mathbb{F}(N)$ is also a projective cover in $\mod H$.
        \item Assume that $\mathbb{F}(f):\mathbb{F}(P) \rightarrowdbl \mathbb{F}(N)$ is a projective cover of $\mathbb{F}(N)$ induced by $P\in\proj\Lambda$ and a homomorphism $f:P\rightarrowdbl N$. It also gives rise to the commutative diagram \rm{(\ref{CD:big_comm_diagram})}.
        \item Under the assumptions in (3), the homomorphism $f:P \rightarrowdbl N$ is also a projective cover in $\CM \Lambda$.
    \end{enumerate}
            \begin{equation} \label{CD:big_comm_diagram}
            \begin{tikzcd}[row sep=1.25em, column sep={between origins, 4em}]
                    &                &                                         &                                           &                                                         &0 \arrow[dd]                                       &                                                      &0 \arrow[dd]                                    &    & \\
                    &                &0 \arrow[dd]                             &                                           &0 \arrow[dd]                                             &                                                   &0 \arrow[dd]                                          &                                                &    & \\
                    &0 \arrow[rr]    &                                         &X \arrow[rr] \arrow[dl,equal] \arrow[dd]   &                                                         &\rad P \arrow[rr] \arrow[dl,equal] \arrow[dd]      &                                                      &\rad N \arrow[rr] \arrow[dl,equal] \arrow[dd]   &    &0\\
    0 \arrow[rr]    &                &X \arrow[rr,crossing over] \arrow{dd}    &                                           &\rad P \arrow[rr,crossing over] \arrow[dd]               &                                                   &\rad N \arrow[rr, crossing over] \arrow[dd]           &                                                &0   & \\
                    &0 \arrow[rr]    &                                         &X \arrow[rr] \arrow[dl]                    &                                                         &P \arrow[rr,"f"near start] \arrow[dl] \arrow[dd]             &                                                      &N \arrow[rr] \arrow[dl] \arrow[dd] \arrow[dd]   &    &0\\
    0 \arrow[rr]    &                &Y \arrow[rr] \arrow[dd]                  &                                           & P \Gamma \arrow[rr,"f_1"near end,crossing over] \arrow[dd]             &                                                   & N \Gamma \arrow[rr,crossing over] \arrow[dd]          &                                                &0   & \\
                    &                &                                         &                                           &                                                         &P/\rad P \arrow[rr,"\bar{f}"near start,"\simeq"'near start] \arrow[dl] \arrow[dd]  &                                            &N/\rad N \arrow[dl]  \arrow[dd]       &    &\\
    0 \arrow[rr]    &                &Z \arrow[rr] \arrow[dd]                  &                                           &P\Gamma/\rad P \arrow[rr,"\bar{f_1}"near start,crossing over] \arrow[dd]      &                                                   &N\Gamma /\rad N \arrow[rr,crossing over] \arrow[dd]  &                                                &0   & \\
                    &                &                                         &                                           &                                                         &0                                                  &                                                      &0                                               &    & \\
                    &                &0                                        &                                           &0                                                        &                                                   &0                                                     &                                                &    & 
            \arrow[from=2-3, to=4-3, crossing over]  
            \arrow[from=2-5, to=4-5, crossing over]  
            \arrow[from=3-4, to=5-4, "\rotatebox{270}{$\simeq$}" near start]
            \arrow[from=4-3, to=6-3, crossing over]  
            \arrow[from=4-3, to=4-5, crossing over]  
            \arrow[from=4-5, to=6-5, crossing over]  
            \arrow[from=4-7, to=6-7, crossing over]  
            \arrow[from=6-7, to=8-7, crossing over]  
        \end{tikzcd}
        \end{equation}
    \begin{proof}
        (1): Assume that $f: P \rightarrowdbl N$ is a projective cover in $\CM \Lambda$. Let $X$ be the kernel of $f$ and $Y$ be the kernel of the induced map $f_1: P \Gamma\rightarrowdbl N\Gamma $. There exists a commutative diagram of exact sequences:
        \begin{equation}\label{CD:bigcommute1}
            \begin{tikzcd}
                0 \arrow[r] &X \arrow[r] \arrow[d,dotted] &P \arrow[r,"f"] \arrow[d] &N \arrow[r] \arrow[d] &0\\
                0 \arrow[r] &Y \arrow[r] &P\Gamma \arrow[r,"f_1"] &N\Gamma \arrow[r] &0
            \end{tikzcd}
        \end{equation}
        Let $Z$ be the kernel of $\bar{f_1}:P\Gamma/\rad P \to N\Gamma/\rad N$. By the 3 $\times$ 3 lemma, there exists a commutative diagram with exact columns and rows:
         \begin{equation}\label{CD:bigcommute2}
            \begin{tikzcd}
                & 0\arrow[d,dotted] & 0\arrow[d] &0 \arrow[d] & \\
                0 \arrow[r] & X \arrow[r] \arrow[d,dotted] &\rad P \arrow[r] \arrow[d] & \rad N \arrow[r] \arrow[d] &0 \\
                0 \arrow[r] & Y \arrow[r] \arrow[d,dotted] &P\Gamma \arrow[r,"f_1"] \arrow[d] & N\Gamma \arrow[r] \arrow[d] &0 \\
                0 \arrow[r] & Z \arrow[r] \arrow[d,dotted] &P\Gamma/\rad P \arrow[r,"\bar{f_1}"] \arrow[d] & N\Gamma/\rad N \arrow[r] \arrow[d] &0 \\
                & 0 &0 &0 & 
            \end{tikzcd}
        \end{equation}
        Combining the diagrams (\ref{cartesian_diagram}), (\ref{CD:bigcommute1}) and (\ref{CD:bigcommute2}), we obtain the desired commutative diagram (\ref{CD:big_comm_diagram}). 

        (2): Note that as an $H$-module, we have $\rad ( P\Gamma  / \rad P, P / \rad P) = (P\Gamma  / \rad P, 0)$. So, the exact sequence in $\mod H$
        \begin{equation*}
            \begin{tikzcd}[column sep=2em]
                0  \arrow[r] & (Z,0) \arrow[r] & ( P \Gamma / \rad P,P / \rad P) \arrow[r,"\mathbb{F}(f)"]     &     ( N \Gamma / \rad N , N / \rad N)       \arrow[r] &    0
            \end{tikzcd}
        \end{equation*}
        is a minimal projective resolution of $\mathbb{F}(M)$. 

        (3)(4): We just need to note that each step above is invertible. In detail, assume that 
        \begin{equation*}
            \mathbb{F}(f): \mathbb{F}(P)=(P\Gamma /\rad P, P/\rad P) \rightarrowdbl (N\Gamma /\rad N, N/\rad N)=\mathbb{F}(N)
        \end{equation*}
        is a projective cover. Then the kernel of $\mathbb{F}(f)$ is of the form $(Z,0)$ and $P/\rad {P} \simeq N/\rad N$. So, again let $X$ be the kernel of $f$ and $Y$ be the kernel of the induced map $ P\Gamma  \rightarrowdbl  N\Gamma $ by $f$, we obtain the desired big commutative diagram {\rm(\ref{CD:big_comm_diagram})}. Finally, since $X=\Ker f \subseteq \rad P$, the homomorphism $f:P \rightarrowdbl N$ is a projective cover in $\CM \Lambda$.
    \end{proof}
\end{proposition}

\begin{corollary}\label{cor:compute syzygy}
    Let $0\to X \to P\xrightarrow{f}M\to0$ be an exact sequence where $f:P\rightarrowdbl M$ is a projective cover for $M\in\CM\Lambda$. Then $\Ker \mathbb{F}(f)$=$((\corad X)/X,0)$ \rm{(recall {Definition} \ref{def:corad} for coradicals)}.
    \begin{proof}
        Consider the commutative diagram (\ref{CD:big_comm_diagram}). The top and middle outer exact sequences split, since they are in $\CM \Gamma$. So, $X\simeq \rad Y$. Thus, $Y \simeq \corad X$ and $Z \simeq Y/X \simeq (\corad X)/X$. 
    \end{proof}
\end{corollary}

The following example is to explain Proposition \ref{thm:syzygy_in_addGamma}

\begin{example}
    Let $(R,\pi R,k)$ be a complete discrete valuation ring.
    \begin{equation*}
    \Lambda=\begin{pNiceMatrix}[last-col] R & R & \Block{2-1}{P_1}\\ \pi R & R & \CodeAfter \tikz \draw ([xshift=0.3em,yshift=0.2em] 1-1.south east) -- ([xshift=-0.1em,yshift=-0.1em] 2-2.north west); \tikz \draw ([xshift=0.2em,yshift=0.05em] 1-1.south east) -- ([xshift=-0.2em, yshift=-0.25em] 2-2.north west);\end{pNiceMatrix} \quad \text{and} \quad \Gamma=\begin{pNiceMatrix}[last-col]  R & R & Q_2\\  \pi R & R & Q_3 \end{pNiceMatrix}      
    \end{equation*} where
        $\left(R=R\right):=\{(r_1,r_2)\in R \times R \mid r_1-r_2 \in \pi R\}$.
    The labels on the right indicate the row vectors, namely, $P_1$ denotes the indecomposable projective $\Lambda$-module, and $Q_2,Q_3$ denote the indecomposable projective $\Gamma$-modules. These constitute all the indecomposable Cohen-Macaulay $\Lambda$-modules up to isomorphism. We have $\corad Q_2 \simeq Q_3$ and $\corad Q_3 \simeq Q_2$. The finite dimensional $k$-algebra $H(\Lambda,\Gamma)$ is isomorphic to the path algebra of the quiver
        $\begin{tikzcd}[column sep=1.5em]
             2 & 1 \arrow[l] \arrow[r] & 3.
        \end{tikzcd}$

    Under the functor $\mathbb{F}$, the indecomposable non-simple projective $H$-module $P(1)$ at 1 and the indecomposable non-simple injective $H$-module $I(2)$, $I(3)$ at 2 and 3 correspond to $P_1$, $Q_2$ and $Q_3$ respectively. Then the minimal projective resolutions in $\mod H$
    \begin{equation*}
        0 \to P(2) \to P(1) \to I(3) \to 0, \quad 0 \to P(3) \to P(1) \to I(2) \to 0,
    \end{equation*}
    give rise to the projective covers in $\CM\Lambda$
    \begin{equation*}
        0 \to Q_2 \to P_1 \to Q_2 \to 0, \quad 0 \to Q_3 \to P_1 \to Q_3 \to 0.
    \end{equation*}
\end{example}

\section{Singularity Categories of B{\"a}ckstr{\"o}m Orders}\label{sec:sing}

We prove our main theorems in this section. Throughout this section, let $\Lambda=(\Lambda,\Gamma)$ be a B{\"a}ckstr{\"o}m order and keep the notations in Section \ref{subsec: The Representation Theory of a B-Order}.

\subsection{Explicit Descriptions of Singularities of B{\"a}ckstr{\"o}m orders}\label{subsec: The Singularity Category of a B-Order}


We first give a concrete description of the syzygy functor. Recall that we have $\Omega^2(\smod \Lambda) \subseteq \Omega(\sCM\Lambda)\subseteq \sadd \Gamma$; see Proposition \ref{prop:OmegaCMLambda}. So, the syzygy functor $\Omega: \smod \Lambda \to \smod \Lambda$ restricts to a functor $\Omega: \sadd \Gamma \to \sadd \Gamma$.  Consider the multiplication map $\mu: \Gamma \otimes_R \Lambda \rightarrowdbl \Gamma$. Since $\Gamma$ is free as an $R$-module, $\Gamma \otimes_R \Lambda$ is projective as a right $\Lambda$-module. Hence, it gives a projective presentation of $\Gamma$ in $\CM \Lambda$:
\begin{equation}\label{eqa: proj. resol.of Gamma}
    0\to\Ker\mu \to \Gamma \otimes_R \Lambda \xrightarrow{\mu} \Gamma \to 0,
\end{equation}
Then $\Ker \mu$ is a $(\Gamma,\Lambda)$-bimodule and $\Omega(\Gamma)=\Ker\mu$. We consider the stable endomorphism ring $\sEnd_\Lambda(\Gamma)$ and the stable Hom-space $\sHom_\Lambda(\Gamma,\Ker\mu)$. 

\begin{lemma}\label{prop:sEndGamma_semisimple}
    The stable endomorphism ring $\sEnd_\Lambda(\Gamma)$ is a finite dimensional semisimple $k$-algebra.
    \begin{proof}
        First, recall that the restriction functor $\res: \CM \Gamma \to \CM \Lambda$ is fully faithful. Since $\rad \Gamma \simeq \Hom_\Gamma(\Gamma,\rad \Gamma)\subseteq \Hom_\Gamma(\Gamma,\Gamma)$, by Lemma \ref{lem:P(M,N)=(M,rad N)}, we have $\rad \Gamma \subseteq \mathcal{P}_\Lambda(\Gamma,\Gamma)$. Thus, the proposition follows from the projection $\End_\Gamma (\Gamma)/\rad \End_\Gamma(\Gamma) \rightarrowdbl \sEnd_\Lambda(\Gamma)$.
    \end{proof}
\end{lemma}

Henceforth, we denote
    \begin{equation*}
        D:=\sEnd_\Lambda(\Gamma) \quad \text{and} \quad M:=\sHom_\Lambda(\Gamma,\Ker\mu),
    \end{equation*}
The following proposition realizes the syzygy functor as a tensor functor.
    
\begin{proposition}\label{prop: M is a (D,D)-bimodule}
       The following statements hold.
        \begin{enumerate}[label=\textnormal{(\arabic*)}]
            \item M is an $(D,D)$-bimodule.
        \item The following diagram commutes up to a natural isomorphism of functors:
        \begin{equation*}
            \begin{CD}
                \sadd \Gamma @>\sHom_\Lambda(\Gamma,?)>\simeq> \mod D  \\
                @VV\Omega V                                     @VV? \otimes_D MV  \\
                \sadd \Gamma @>\sHom_\Lambda(\Gamma,?)>\simeq> \mod D 
            \end{CD}
        \end{equation*}
        \end{enumerate}
        \begin{proof}
            The Hom-space $M=\sHom_\Lambda(\Gamma, \Ker\mu)$ is a right $\sEnd_\Lambda (\Gamma)$-module via composition of maps. In addition, the syzygy functor gives a map $\Omega: \sEnd_\Lambda (\Gamma) \to \sEnd_\Lambda(\Omega(\Gamma))= \sEnd_\Lambda (\Ker\mu)$, which induces a left $\sEnd_\Lambda (\Gamma)$-module structure on $M$. These indeed define a $(D,D)$-bimodule structure on $M$. 
            
            For (2), the vertical equivalence is by projectivization and $\proj D=\mod D$, since $D$ is semisimple by Lemma \ref{prop:sEndGamma_semisimple}. For the commutative diagram, we first construct a natural transformation between these two functors. Let $Q\in\CM\Gamma$. Tensoring $Q$ with the exact sequence (\ref{eqa: proj. resol.of Gamma}), we obtain an exact sequence
            \begin{equation*}
                0 \to Q\otimes_\Gamma \Ker \mu \to Q\otimes_R \Lambda \to Q \to 0.
            \end{equation*}
            Thus, $\Omega(Q)=Q\otimes_\Gamma \Ker\mu$. We then have a natural isomorphism 
            \begin{equation*}
                \Hom_\Lambda(\Gamma,Q)\otimes_\Gamma M \simeq Q\otimes_\Gamma M \xrightarrow{\simeq} \sHom_\Lambda(\Gamma,Q\otimes_\Gamma\Ker\mu),
            \end{equation*}
            where the first isomorphism is by Proposition \ref{prop:overorder} and the second isomorphism is because $Q\in\proj\Gamma$. On the other hand, we have a natural projection 
            \begin{equation*}
                \Hom_\Lambda(\Gamma,Q)\otimes_\Gamma M \rightarrowdbl \sHom_\Lambda(\Gamma,Q)\otimes_D M,
            \end{equation*}
            which is induced by $\Gamma \xrightarrow{\simeq} \End_\Gamma(\Gamma) \rightarrowdbl D$. Combining together with the natural isomorphism above, we obtain the desired natural transformation. The natural transformation is an isomorphism when evaluated at $Q=\Gamma$. Therefore, we obtain (2). 
        \end{proof}
    \end{proposition}

\begin{proposition}\label{eqa:left triangulated equivalence}
    The categories $(\sadd \Gamma, \Omega)$ and $(\mod D, ?\otimes_D M)$ are trivial left triangulated categories. The functor $\sHom_\Lambda(\Gamma,?):(\sadd \Gamma, \Omega)\to (\mod D, ?\otimes_D M)$ gives a left triangle equivalence.
    \begin{proof}
        Since $\Omega^2(\smod\Lambda)\subseteq \sadd\Gamma$, it is easy to see that $(\sadd \Gamma, \Omega)$ forms a left triangulated subcategory of $(\smod \Lambda,\Omega)$. By Lemmas \ref{prop:sEndGamma_semisimple} and \ref{lem:semisimpletriangle}, the left triangulated structure of $(\sadd \Gamma, \Omega)$ is trivial. By the same lemmas, $(\mod D, ?\otimes_D M)$ is also a trivial left triangulated category. Therefore, Proposition \ref{prop: M is a (D,D)-bimodule} implies that $\sHom_\Lambda(\Gamma,?)$ is a left triangle equivalence.
    \end{proof}
\end{proposition}

We are now prepared to prove the main theorem in this article. Recall Definition-Theorem \ref{ex:V(D,M)}, define
\begin{equation*}
    V(\Lambda):=V(D,M)=\underset{i\geq0}{\varinjlim} \End_D(M^{\otimes i}) \quad {\rm{and}} \quad K(\Lambda):=K(D,M)=\underset{i\geq1}{\varinjlim} \Hom_D(M^{\otimes i}, M^{\otimes i-1}).
\end{equation*}
Then $V(\Lambda)$ is a von Neumann regular algebra. The category $\proj V(\Lambda)$ of finitely generated projective $V(\Lambda)$-modules is semisimple abelian, and $(\proj V(\Lambda),? \otimes_{V(\Lambda)} K(\Lambda))$ forms a trivial left triangulated category.

\begin{theorem}\label{thm:DsgBack}
    For a B{\"a}ckstr{\"o}m order $(\Lambda,\Gamma)$, there exists a triangle equivalence
    \begin{equation*}
        (\Dsg(\Lambda),[1]) \simeq (\proj V(\Lambda), ? \otimes_{V(\Lambda)}  K(\Lambda)),
    \end{equation*}
    where the triangulated structures are trivial.
    \begin{proof}
        By Propositions \ref{prop:OmegaCMLambda} and \ref{prop:Dsg_iso_projEnd}, we have $\Dsg(\Lambda) \simeq \proj \End_{\Dsg(\Lambda)}(q(\Gamma))$ as additive categories. On the other hand, Theorem \ref{thm:stabilization} and Proposition \ref{prop: M is a (D,D)-bimodule} imply 
        \begin{equation*}
            \End_{\Dsg(\Lambda)}(q(\Gamma)) \simeq \underset{i\geq0}{\varinjlim} \sEnd_\Lambda(\Omega^i (\Gamma))
                                          \simeq \underset{i\geq0}{\varinjlim} \End_D(M^{\otimes i})=V(\Lambda).
        \end{equation*}
        Therefore, $\Dsg(\Lambda) \simeq \proj V(\Lambda)$ as additive categories. Let $\Sigma$ be the corresponding automorphism to $[1]$. Then $(\proj V(\Lambda),\Sigma)$ becomes a triangulated category inherited by the categorical equivalence. In addition, by Lemma \ref{lem:semisimpletriangle}, there is a unique triangulated structure with respect to $\Sigma$, which is the trivial triangulated structure. So, these two triangulated structures coincide. Therefore, we obtain a triangle equivalence
        \begin{equation*}
            \Hom_{\Dsg(\Lambda)}(q(\Gamma),?):(\Dsg(\Lambda),[1]) \xlongrightarrow{\sim} (\proj V(\Lambda), \Sigma).
        \end{equation*}

        Next, we determine the suspension functor $\Sigma$. Note that $\Sigma(V(\Lambda))$ is a $(V(\Lambda),V(\Lambda))$-bimodule whose left module structure is given by $V(\Lambda) \simeq \End_{V(\Lambda)}(V(\Lambda)) \xrightarrow{\Sigma} \End_{V(\Lambda)}(\Sigma(V(\Lambda)))$. Moreover, we have a natural isomorphism $\Sigma \simeq ? \otimes_{V(\Lambda)} \Sigma (V(\Lambda)): \proj V(\Lambda) \to \proj V(\Lambda)$.

        A similar calculation shows that there exists a $(V(\Lambda),V(\Lambda))$-bimodule isomorphism
        \begin{equation*}
            \Sigma(V(\Lambda))\simeq\Hom_{\Dsg(\Lambda)}(q(\Gamma),q(\Gamma) [1]) \simeq \underset{i\geq1}{\varinjlim} \sHom_\Lambda(\Omega^i(\Gamma), \Omega^{i-1} (\Gamma))
                                               \simeq \underset{i\geq1}{\varinjlim} \Hom_D(M^{\otimes i}, M^{\otimes i-1})=K(\Lambda).
        \end{equation*}
       Therefore, $\Sigma$ is isomorphic to the functor $? \otimes_{V(\Lambda)} K (\Lambda): \proj V(\Lambda) \to \proj V(\Lambda)$. 
    \end{proof}
\end{theorem}

As noted in the introduction, our result on the singularity category of a B{\"a}ckstr{\"o}m order closely parallels Chen's result; see Theorem \ref{thm:DsgArtin}. Using the notations above, we associate to a B{\"a}ckstr{\"o}m order $\Lambda=(\Lambda,\Gamma)$ a finite dimensional $k$-algebra $A(\Lambda)$ defined as the trivial extension 
\begin{equation}\label{eqa:A(Lambda)}
    A(\Lambda):=D\oplus M,
\end{equation}
where the multiplication is given by $(a,m)(a',m'):=(aa',am'+ma')$. Since $D$ is semisimple and $M^2=0$ in $A(\Lambda)$, we have $\rad A(\Lambda)\simeq M$ and $A(\Lambda)/\rad A(\Lambda)\simeq D$. Therefore, $A(\Lambda)$ is a finite dimensional radical square zero $k$-algebra.

\begin{corollary}\label{cor:DsgLambda=DsgA}
    Let $\Lambda$ be a B{\"a}ckstr{\"o}m order and $A(\Lambda)=D\oplus M$ be the finite dimensional radical square zero algebra associated to it. We have a triangle equivalence 
    \begin{equation*}
        \Dsg(\Lambda) \simeq \Dsg(A(\Lambda)).
    \end{equation*}
    \begin{proof}
        By construction, we have $V(\Lambda)\simeq V(A(\Lambda))$ and $K(\Lambda) \simeq K(A(\Lambda))$. Hence, Theorems \ref{thm:DsgBack} and \ref{thm:DsgArtin} yield the desired triangle equivalence. 
    \end{proof}
\end{corollary}

We also construct the equivalence functor explicitly using the theory of stabilizations. Recall that $\sadd D=\sadd_{A(\Lambda)}D$ denotes the full subcategory of $\smod A(\Lambda)$ generated by $D$.

\begin{proposition}\label{prop:explicit equivalence functor}
    Let $\Lambda$ be a B{\"a}ckstr{\"o}m order and $A:=A(\Lambda)=D\oplus M$ be the finite dimensional radical square zero algebra associated to it. Then 
    \begin{enumerate}[label=\textnormal{(\arabic*)}]
        \item There is a natural full and dense left triangle functor between two trvial left triangulated categories
        \begin{equation*}\pi: (\mod D, ? \otimes_D M) \longrightarrow (\sadd D, ? \otimes_D M).
        \end{equation*}
        \item The functor $\pi$ induces a triangle equivalence $\tilde{\pi}: \mathcal{S}(\mod D, ? \otimes_D M) \xlongrightarrow{\simeq} \mathcal{S}(\sadd D, ? \otimes_D M)$ such that the following diagram commutes:
    \begin{equation}\label{CD: explicit equivalence functor}
        \begin{CD}
            (\mod D, ? \otimes_D M) @>S>> \mathcal{S}(\mod D,? \otimes_D M)\\
            @VV\pi V @V\rotatebox{90}{\tiny $\simeq$} V\tilde{\pi}V \\
            (\sadd D, ? \otimes_D M) @>S>> \mathcal{S}(\sadd D,? \otimes_D M)
        \end{CD}
    \end{equation}
        \item There are triangle equivalences $\Dsg(\Lambda)\simeq\mathcal{S}(\mod D, ?\otimes_D M)$ and $\Dsg(A)\simeq\mathcal{S}(\sadd D, ?\otimes_D M)$. Hence, the functor $\tilde{\pi}$ gives the equivalence $\Dsg(\Lambda)\xlongrightarrow{\simeq} \Dsg(A)$.
    \end{enumerate}
    \begin{proof}
        (3): By Theorem \ref{thm:stabilization}, we have $\Dsg(\Lambda) \simeq \mathcal{S}(\smod \Lambda, \Omega_\Lambda)$ and $\Dsg(A) \simeq \mathcal{S}(\smod A, \Omega_A)$. Then, by Lemma \ref{pro:criterion for tilde{F}}, we have $\mathcal{S}(\smod \Lambda, \Omega_\Lambda)\simeq \mathcal{S}(\sadd \Gamma, \Omega_\Lambda)$ and $\mathcal{S}(\smod A, \Omega_A) \simeq \mathcal{S}(\smod A, \Omega_A)$. Finally, by Proposition \ref{eqa:left triangulated equivalence} and Theorem \ref{thm:DsgArtin}, we have $\mathcal{S}(\sadd \Gamma, \Omega_\Lambda) \simeq \mathcal{S}(\mod D, ?\otimes_D M)$ and $\mathcal{S}(\smod A, \Omega_A)\simeq \mathcal{S}(\sadd D, ?\otimes_D M)$.
        
        (1): Since $D \simeq A/\rad A$, we have $\mod D=\add_D D=\add_A D$. Thus there is a canonical projection $\pi: \mod D \to \sadd_A D$ which is full and dense. By Theorems \ref{thm:DsgBack} and \ref{thm:DsgArtin}, the left triangle structures of $(\mod D, ?\otimes_D M)$ and $(\sadd_A D, ? \otimes_D M)$ are trivial. Therefore, $\pi$ is a left triangle functor.
        
        (2): By the universal property of stabilization, $\pi$ induces a triangle functor $\tilde{\pi}$ between their stabilizations such that the commutative diagram (\ref{CD: explicit equivalence functor}) holds. Now, we show that $\tilde{\pi}$ is an equivalence. Since $\pi$ is full and dense, by Lemma \ref{pro:criterion for tilde{F}} (1)(3), $\tilde{\pi}$ is also full and dense. For the faithfulness, we consider two morphisms $f,f': S \to T$ between two simple $A$-modules $S,T$ such that $\pi(f)=\pi(f')$. If $S \not\simeq T$, since $\Hom_A(S,T)=0$, trivially $f=f'=0$. If $S\simeq T \not\in \proj A$, then $\mathcal{P}(S,S)=0$ and hence $\Hom_A(S,S)=\sHom_A(S,S)$. Thus, $f=\pi(f)=\pi(f')=f'$. If $S\simeq T \in \proj A$, then $S\otimes_D M \simeq \Omega(S) \simeq \rad S$=0. Therefore, $f \otimes_D M = f' \otimes_D M=0$. By Lemma \ref{pro:criterion for tilde{F}} (2), $\tilde{{\pi}}$ is faithful. 
    \end{proof}
\end{proposition}

\subsection{Classifications of some Classes of B{\"a}ckstr{\"o}m Orders}\label{subsec: Classifications of some Classes of B-Orders}

We apply Corollary \ref{cor:DsgLambda=DsgA} to give criteria for weakly regular, Gorenstein, Iwanaga-Gorenstein and sg-Hom-finite B{\"a}ckstr{\"o}m orders via the associated finite dimensional radical square zero algebras. The relationships between these classes of B{\"a}ckstr{\"o}m orders are as presented in the introduction; see the figure (\ref{eqa:hierarchy}).

\begin{definition}\cite[Page 69]{ARS97}.
    Let $A$ be a basic finite dimensional algebra. Assume that $A/\rad A \simeq \prod_{i\in I}D_i$, where $I$ is a finite index set and $D_i$ are finite dimensional division algebras. Let $\{e_i\}_{i\in I}$ be a complete set of primitive orthogonal idempotents of $A$ in which $e_i$ corresponds to $D_i$. Let $M_{ij}:=e_i(\rad A/\rad^2 A)e_j$ be a $(D_i,D_j)$-bimodule for any $i,j\in I$. The \textit{valued quiver} $Q_A$ is defined as follows. 
    
    The vertex set is $I$, corresponding to the division algebras $\{D_i\}_{i\in I}$. For any $i,j\in I$, there is an arrow from $i$ to $j$ endowed with a \textit{valuation} $(\dim (M_{ij})_{D_j}, \dim _{D_i}(M_{ij}))$.  For a general finite dimensional algebra $A$, its valued quiver $Q_A$ is defined as the valued quiver of its corresponding basic algebra.
    
    A valuation is called \textit{trivial} if it is (1,1). The valued quiver $Q_A$ is called \textit{acyclic} if it contains no oriented cycle.

\end{definition}

We first recall the results by Xiao-wu Chen on the classifications of the similar classes of a finite dimensional radical square zero algebras in terms of their valued quivers.

\begin{proposition} \cite[Theorem 5.2]{Che11} \cite[Corollary 1.3]{Che12}  \label{prop:criterion for IG artinian} \label{thm:HomfiniteDsgArtin}
    Let $A$ be a finite dimensional radical square zero algebra and $Q_A$ be its valued quiver. The following statements hold.
    \begin{enumerate}[label=\textnormal{(\arabic*)}]
        \item $A$ is weakly regular if and only if $Q_A$ is acyclic.
        \item $A$ is self-injective if and only each connected component of $Q_A$ is a single vertex or a cycle with trivial valuations.
        \item $A$ is Iwanaga-Gorenstein if and only if each connected component of $Q_A$ is acyclic or is a cycle with trivial valuations.
        \item $A$ is sg-Hom-finite if and only if $V(A)$ is semisimple if and only if $Q_A$ is obtained from a disjoint union of oriented cycles with trivial valuations by adjoining sources and sinks with arbitrary values repeatedly.
    \end{enumerate}
\end{proposition}

Let $\Lambda$ be a B{\"a}ckstr{\"o}m order and $A(\Lambda)$ be the finite dimensional radical square zero $k$-algebra associated to it; recall (\ref{eqa:A(Lambda)}). 

\begin{corollary}\label{cor:HomfiniteDsgBack}
    The following statements are equivalent.
    \begin{enumerate}[label=\textnormal{(\arabic*)}]
        \item $\Lambda$ is sg-Hom-finite;
        \item $A(\Lambda)$ is sg-Hom-finite;
        \item $V(\Lambda)=V(A(\Lambda))$ is semisimple.
        \item $Q_{A(\Lambda)}$ is obtained from a disjoint union of oriented cycles with trivial valuations by adjoining sources and sinks with arbitrary values repeatedly.
    \end{enumerate}
    \begin{proof}
       $(1) \Leftrightarrow (2)$ is by Corollary \ref{cor:DsgLambda=DsgA}. $(2) \Leftrightarrow (3) \Leftrightarrow (4)$ is by Proposition \ref{prop:criterion for IG artinian} (4).
    \end{proof}
\end{corollary}

\begin{corollary}\label{cor:glfiniteBack}
    The following statements are equivalent.
    \begin{enumerate}[label=\textnormal{(\arabic*)}]
        \item $\Lambda$ is weakly regular;
        \item $\Dsg(\Lambda)=0$;
        \item $A(\Lambda)$ is weakly regular.
        \item $\Dsg(A(\Lambda))=0$;
        \item $Q_{A(\Lambda)}$ is acyclic.
    \end{enumerate}
    \begin{proof}
        $(2) \Leftrightarrow (4)$ is by Corollary \ref{cor:DsgLambda=DsgA}. $(1) \Leftrightarrow (2)$ and $(3) \Leftrightarrow (4)$ is by the fact that the singularity category vanishes if and only if the ring is of finite global dimension. $(3) \Leftrightarrow (5)$ is by Proposition \ref{prop:criterion for IG artinian} (1).
    \end{proof}
\end{corollary}

Next, we classify Iwanaga-Gorenstein B{\"a}ckstr{\"o}m orders. The following property is basic.

\begin{lemma}\cite[Lemma 2.2]{Che12} \label{lem:OmegaofindGp}
    Let $\Lambda$ be any $R$-orders. Let $Q$ be a indecomposable non-projective Gorenstein-projective $\Lambda$-module. Then the kernel of projective cover of $Q$ is also indecomposable non-projective Gorenstein-projective.
\end{lemma}

\begin{proposition}\label{prop:restriction functor to GP}
    Assume that $\Lambda$ is a B{\"a}ckstr{\"o}m order and $A:=A(\Lambda)=D\oplus M$ is the finite dimensional radical square algebra associated to it. Let $\pi,\tilde{\pi}$ be the functors in the commutative diagram {\rm{(\ref{CD: explicit equivalence functor})}} of {\rm{Proposition \ref{prop:explicit equivalence functor}}}. Then $\pi$ restricts to a triangle equivalence $\pi':\sGp \Lambda \to \sGp A$. More precisely, we have the following commutative diagram:
    \begin{equation*}
            \begin{tikzcd}
                (\sGp \Lambda, \Omega_\Lambda) \arrow[d,"\pi'","\rotatebox{90}{\tiny $\simeq$}"'] \arrow[r,hookrightarrow, "\subseteq"] & (\sadd_\Lambda \Gamma, \Omega_\Lambda) \arrow[r, "\simeq"] & (\mod D, ? \otimes_D M) \arrow[r, "S"] \arrow[d,"\pi"] & \mathcal{S}(\mod D,? \otimes_D M) \arrow[d,"\tilde{\pi}","\rotatebox{90}{\tiny $\simeq$}"']\\
                (\sGp A,\Omega_A) \arrow[r, hookrightarrow, "\subseteq"] & (\sadd_A D, \Omega_A) \arrow[r,"\simeq"]& (\sadd D, ? \otimes_D M) \arrow[r,"S"] & \mathcal{S}(\sadd D,? \otimes_D M)
            \end{tikzcd}
    \end{equation*}
    \begin{proof}
        Recall that by Proposition \ref{eqa:left triangulated equivalence} and Theorem \ref{thm:DsgArtin}, we have left triangle equivalences $a: (\sadd_\Lambda \Gamma, \Omega_\Lambda) \xrightarrow{\simeq}(\mod D, ? \otimes_D M)$ and $b: (\sadd_A D, \Omega_A) \xrightarrow{\simeq} (\sadd D, ?\otimes_D M)$. Recall that $\sGp \Lambda$ and $\sGp A$ are triangulated subcategories of $(\sadd_\Lambda \Gamma, \Omega_\Lambda)$ and $(\sadd_A D, \Omega_A)$ respectively; see Section \ref{subsec:Gorenstein-projective modules}.

        Firstly, we show that the essential image $b^{-1} \pi a (\sGp \Lambda)$ is contained in $\sGp A$, and thus the triangle functor $\pi':=b^{-1}\pi a : \sGp \Lambda \to \sGp A$ is well-defined. Since $\pi$ is a full left triangle functor, $b^{-1} \pi a (\sGp \Lambda)$ is a triangle subcategory of $(\sadd_A D, \Omega_A)$ and hence of $(\smod A, \Omega_A)$. By Proposition \ref{prop: maximum property of GP}, we have $b^{-1} \pi a (\sGp \Lambda) \subseteq \sGp A$. Thus, we have the following commutative diagram.
        \begin{equation*}
            \begin{tikzcd}
                (\sGp \Lambda, \Omega_\Lambda) \arrow[r,hookrightarrow,"\subseteq"] \arrow[d,dotted, "\pi'"]& (\sadd_\Lambda \Gamma, \Omega_\Lambda) \arrow[r, "a", "\simeq"'] & (\mod D, ? \otimes_D M) \arrow[r, "S"] \arrow[d, "\pi"] & \mathcal{S}(\mod D,? \otimes_D M) \arrow[d,"\tilde{\pi}","\rotatebox{90}{\tiny $\simeq$}"']\\
                (\sGp A,\Omega_A) \arrow[r, hookrightarrow, "\subseteq"] & (\sadd_A D, \Omega_A) \arrow[r, "b", "\simeq"']& (\sadd D, ? \otimes_D M) \arrow[r,"S"] & \mathcal{S}(\sadd D,? \otimes_D M)
            \end{tikzcd}
        \end{equation*}

        Secondly, we show that $\pi'$ is full and faithful. Trivially, $\pi'$ is full since $\pi$ is full. To prove the faithfulness, we consider a morphism $f: X \to Y$ between two indecomposable non-projective Gorenstein-projective $\Lambda$-modules $X$ and $Y$ such that $b^{-1}\pi a (f)=0$. We will show that $a(f)=0$. Note that $a(X)$ and $a(Y)$ are simple $A$-modules. If $a(X) \not\simeq a(Y)$, then clearly $a(f)=0$. If $a(X)\simeq a(Y) \not\in \proj A$, then $\Hom_A(a(X),a(Y))=\sHom_A(a(X),a(Y))$, and thus $a(f)=\pi a(f)=0$. If $a(X)\simeq a(Y) \in \proj A$, then $a(\Omega_\Lambda(X)) \simeq a(X)\otimes_D M =\rad a(X)=0$, which is a contradiction to $\pd X=\infty$.

        Thirdly, we show that $\pi'$ is dense. Consider the decomposition of $D$ as an $A$-module $D=D_0 \oplus D_1$ such that $D_1 \in \proj A$ and each direct summand of $D_0$ is not $A$-projective. Then there is an equivalence $c: \add D_0 \xrightarrow{\simeq} \sadd D$ as additive categories. Moreover, the inclusion functor $i: \add D_0 \hookrightarrow \mod D$ is left adjoint to the projection $c^{-1}\pi: \mod D \to \add D_0$. By Lemma \ref{lem:OmegaofindGp}, the essential image $ic^{-1} b (\sGp A)$ is a triangle subcategory of $(\mod D, \otimes_D M)$. So, $a^{-1}ic^{-1} b (\sGp A)$ is a triangle subcategory of $(\sadd_\Lambda \Gamma, \Omega_\Lambda)$ and hence of $(\smod \Lambda, \Omega_\Lambda)$. By Proposition \ref{prop: maximum property of GP}, we have $a^{-1}ic^{-1} b (\sGp A) \subseteq \sGp \Lambda$. Then for any $X\in \sGp A$, we have
        \begin{equation*}
            \pi'a^{-1}ic^{-1}ib(X)=(b^{-1}\pi a)(a^{-1}ic^{-1}ib)(X)\simeq X,
        \end{equation*}
        which implies that $\pi'$ is dense.
    \end{proof}
\end{proposition}

\begin{remark}
    In the proof, $\add D_0 \subseteq \mod D$ is not necessarily a left triangulated subcategory. Because even for a non-projective simple $A$-module $S$, $\Omega(S)$ may have an $A$-projective direct summand. However, by Lemma \ref{lem:OmegaofindGp}, this is true for $ic^{-1} b (\Gp A) \subseteq \mod D$.
\end{remark}

\begin{theorem}\label{cor:restriction to GP}
    Assume that $\Lambda$ is a B{\"a}ckstr{\"o}m order and $A:=A(\Lambda)=D\oplus M$ is the finite dimensional radical square algebra associated to it.  The following statements are equivalent.
    \begin{enumerate}[label=\textnormal{(\arabic*)}]
        \item  $\Lambda$ is Iwanaga-Gorenstein;
        \item  $A(\Lambda)$ is Iwanaga-Gorenstein;
        \item  Each connected component in the species of $Q_{A(\Lambda)}$ is acyclic or is a cycle with trivial valuations.
    \end{enumerate}
    \begin{proof}
        By Proposition \ref{prop:restriction functor to GP}, the canonical functor $\sGp \Lambda \to \Dsg (\Lambda)$ is a triangle equivalence if and only if the canonical functor $\sGp A \to \Dsg (A)$ is a triangle equivalence. Then (1) $\Leftrightarrow$ (2) follows by Theorem \ref{thm:criterion for Gor order}. (2) $\Leftrightarrow$ (3) is by Proposition \ref{prop:criterion for IG artinian} (3).
    \end{proof}
\end{theorem}

We finally classify Gorenstein B{\"a}ckstr{\"o}m orders. 

\begin{theorem}\label{cor:criterion for Gor}
    Assume that $\Lambda$ is a B{\"a}ckstr{\"o}m order and $A(\Lambda)=D\oplus M$ is the finite dimensional radical square zero algebra associated to it. The following statements are equivalent.
    \begin{enumerate}[label=\textnormal{(\arabic*)}]
        \item $\Lambda$ is Gorenstein;
        \item $A(\Lambda)$ is a product of non-simple self-injective algebras;
        \item Each connected component of $Q_{A(\Lambda)}$ is a cycle with trivial valuations.
    \end{enumerate}
    
    \begin{proof}
        (1) $\Rightarrow$ (2): Let $\Lambda$ be Gorenstein. Then we have $\CM \Lambda=\Gp \Lambda$; see \cite[Corollary 11.5.3]{EJ00}. So, each non-projective Cohen-Macaulay module has an infinite projective dimension, which implies that each component of $Q_{A(\Lambda)}$ cannot be acyclic. Note that $\Lambda$ is in particular Iwanaga-Gorenstein. Then, by Theorem \ref{cor:restriction to GP}, each connected component of $Q_{A(\Lambda)}$ is a cycle with trivial valuations. So, $A(\Lambda)$ is a product of some non-simple self-injective algebras.

        (2) $\Rightarrow$ (1): Assume conversely that $A(\Lambda)$ is a product of some non-simple self-injective (i.e. non-simple Gorenstein) algebras. Then $\sGp A(\Lambda) \subseteq \sadd D \subseteq \smod A(\Lambda) $ is an equality. To show that $\Lambda$ is Gorenstein, it is equivalently to show $\proj\Lambda = \inj \Lambda$; see \cite [Lemma 2.15]{IW 14}. Without loss of generality, we may assume $\Lambda$ to be ring-indecomposable and non-hereditary. Since each projective $A(\Lambda)$-module is not simple, we have $\sadd D = \mod D$. Therefore, $\sGp \Lambda = \sadd \Gamma$ by Proposition \ref{prop:restriction functor to GP}.
        
        We claim that $\proj \Lambda \cap \proj \Gamma=\{0\}$. In fact, if there exists $0\neq Q\in\ind(\proj\Lambda \cap \proj\Gamma)$, then by the assumptions that $\Lambda$ is ring-indecomposable and $\Lambda \neq \Gamma$, we may assume that $\corad Q\in\ind(\proj\Gamma)\backslash\proj\Lambda$. Otherwise, using the notions in \cite[Theorem 1.6]{HN94}, the unique $\Gamma$-composition series $Q\subsetneq \corad Q \subsetneq \corad^2(Q) \subsetneq \dots \subsetneq \corad^i(Q)\simeq Q$ with $\corad^j (Q)\in\ind(\proj\Lambda \cap \proj\Gamma)$, is also the unique $\Lambda$-composition series of $Q$. Then $\Lambda=\Gamma$, which is a contradiction.

        Since $\corad Q \not \in \proj\Lambda$, $\mathbb{F}(\corad Q)$ is a non-projective injective $H$-module (recall Theorem \ref{thm:rep_of_Back}). Since $H$ is a finite dimensional hereditary radical square zero $k$-algebra and $\Lambda/\rad \Lambda \neq \Gamma/\rad \Gamma$, there exists a projective cover $p: P \to I$ in $\mod H$ of some indecomposable non-projective injective $H$-module $I$ such that $\soc \mathbb{F}(\corad Q)=(\corad Q/Q,0)$ is a direct summand of the kernel of $p: P\rightarrowdbl I$, namely, we have an exact sequence in $\mod H$:
        \begin{equation*}
            0 \to (\corad Q/Q,0) \oplus S \to P \xrightarrow{p} I \to 0,
        \end{equation*}
        for some semisimple projective $H$-module $S$.
        Then Corollary \ref{cor:compute syzygy} implies that there exists an exact sequence of $\Lambda$-modules:
        \begin{equation*}
            0 \to Q \oplus \tilde{S} \to \tilde{P} \xrightarrow{\tilde{p}} \tilde{I} \to 0,
        \end{equation*}
        where $\mathbb{F}(\tilde{p})=p$ is a projective cover of some $\tilde{I} \in \ind(\proj\Gamma)$. Since $\sGp \Lambda = \sadd \Gamma$, $\tilde{I}\in \ind(\Gp \Lambda)$. Since $Q\in\proj\Lambda$, then $Q=0$ by Lemma \ref{lem:OmegaofindGp}, which is a contradiction. This finishes the proof of our claim.

        Now, Proposition \ref{prop:OmegaCMLambda} and our claim imply that every non-projective $\Lambda$-module has infinite projective dimension. Let $N \in \inj \Lambda$. By Lemma \ref{lem:lem to Gor}, we have $\Gpd N=\pd N$. If $N$ is not projective, then $\Gpd N=\infty$, which leads to a contradiction to $\sGp \Lambda = \sadd \Gamma$. Therefore, $\inj\Lambda\subseteq\proj\Lambda$. Finally, since we have the same number of the isomorphism classes of indecomposable projective and injective objects, we have $\inj\Lambda=\proj\Lambda$.
    \end{proof}
\end{theorem}

\section{Examples}\label{sec: example}

Throughout this section, we assume that $(R,\mathfrak{m},k)$ is a complete discrete valuation ring and $\mathfrak{m}=\pi R$. Let $\Lambda=(\Lambda,\Gamma)$ be a B{\"a}ckstr{\"o}m order and keep the notations in Sections \ref{subsec: The Representation Theory of a B-Order} and \ref{sec:sing}. In this section, we present several examples to illustrate our results. 

To compute the associated finite dimensional radical square zero $k$-algebra $A(\Lambda)$, we need to compute the semisimple algebra $D=\sEnd_\Lambda(\Gamma)$ and the $(D,D)$-bimodule $M=\sHom_\Lambda(\Gamma, \Ker \mu)$, where $\mu: \Gamma \otimes_R \Lambda \rightarrowdbl \Gamma$ is the multiplication map. 

Let $\Omega(\Gamma)$ be the kernel of the projective cover of $\Gamma$ in $\CM \Lambda$. Then, as a right $D$-module, we have
\begin{equation*}
    M\simeq \sHom_\Lambda(\Gamma,\Omega(\Gamma)),
\end{equation*}
which can be computed via Proposition \ref{prop:sEndGamma_semisimple}. The left $D$-module structure is given by the $k$-algebra homomorphism on the endomorphism rings induced by the syzygy functor 
\begin{equation*}
    \Omega: D=\sEnd_\Lambda(\Gamma) \to \sEnd_\Lambda(\Omega(\Gamma)).
\end{equation*}

We start with an easy example, Example \ref{ex:1}. An example involving a field extension is given in Example \ref{ex:2}. Example \ref{ex:3} illustrates the hierarchy discussed in the introduction.

\begin{example}\label{ex:1}
    Let $\Gamma=k[\![x_1]\!] \times k[\![x_2]\!] \times \dots \times k[\![x_{n}]\!]$ and $x=(x_1, x_2, \dots, x_{n})$ for $n\geq 2$. Let $\Lambda=k+x\Gamma$. Denote $P_{n+1}=\Lambda$, $Q_{j}=k[\![x_j]\!]$ and $k_j=Q_j/\rad Q_j\simeq k$ for $1\leq j\leq n$. The projective cover of $Q_{j}$ is 
    \begin{equation*}
        0 \to \bigoplus_{j'\neq j} Q_{j'} \to P_{n+1} \to Q_j \to 0.
    \end{equation*}
    Thus, as a right $D$-module, $M=\bigoplus^n_{j=1}(k^{\oplus n-1}_j)$. The syzygy functor is the diagonal map 
    \begin{equation*}
        \begin{CD}
            \Omega: \quad @. \sEnd_\Lambda (Q_j) @>>> \sEnd_\Lambda (\Omega(Q_j))\\
            @.         @|                           @|                              \\
            \Delta: \quad        @.  k_j                     @>>>  \prod_{j'\neq j}k_{j'}        
        \end{CD}
    \end{equation*}
    Represented as matrices, the semisimple algebra $D$ and the $(D,D)$-bimodule $M$ are given by
    \begin{equation*}
        D=
        \begin{pmatrix}
            k & 0 & \cdots & 0& 0 \\
            0 & k & \cdots & 0& 0 \\
            \vdots& \vdots & \ddots & \vdots & \vdots\\
            0 & 0 & \cdots & k& 0 \\
            0 & 0 & \cdots & 0& k
        \end{pmatrix},
        \quad 
        M =
        \begin{pmatrix}
            0 & k & \cdots &k &k \\
            k & 0 & \cdots &k &k\\
            \vdots&\vdots & \ddots & \vdots & \vdots\\
            k & k & \cdots &0 &k\\
            k & k & \cdots &k &0
        \end{pmatrix}.
    \end{equation*}
    
    Therefore, $Q_{A(\Lambda)}$ is a quiver with vertices $1\leq j \leq n$ and exactly one arrow in each direction between any two distinct vertices. The finite dimensional $k$-algebra with radical square zero $A(\Lambda)$ is the path algebra $kQ_{A(\Lambda)}$ with relations that all paths with length at least two vanish. For instance, when $n=3$, then 
    \begin{equation*}
        Q_{A(\Lambda)}:
        \begin{tikzcd}[column sep=1em]
            & 1 \arrow[ld, shift left] \arrow[rd,shift left] & \\
            2 \arrow[ru] \arrow[rr,shift left]& & 3 \arrow[lu] \arrow[ll]
        \end{tikzcd}
    \end{equation*}

    In this example, $\Lambda$ is sg-Hom-finite if and only if $\Lambda$ is Iwanaga-Gorenstein if and only it is Gorenstein if and only if $n=2$.
\end{example}

\begin{example}\label{ex:2}
    Assume $l$ to be a field extension of $k$ with degree $n\geq 2$. Put 
    \begin{equation*}
        \Lambda=k+xl[\![x]\!],\qquad \Gamma=l[\![x]\!],
    \end{equation*}
    According to Proposition \ref{prop:sEndGamma_semisimple}, the semisimple $k$-algebra $D=\sEnd_\Lambda(\Gamma)=\Gamma/\rad\Gamma=l$.
    Consider the exact sequence (\ref{eqa: proj. resol.of Gamma})
    \begin{equation*}
        0 \to \Ker\mu \to \Gamma \otimes_R \Lambda \xrightarrow{\mu} \Gamma \to 0.
    \end{equation*}
    \begin{claim}
        Consider an exact sequence
        \begin{equation*}
            0 \to Z \to l \otimes_k l \xrightarrow{\mu'} l \to 0,
        \end{equation*}
        where $\mu'$ is the multiplication map. Then the $(D,D)$-bimodule $M=\sHom_\Lambda(\Gamma,K)$ is isomorphic to $Z$, where $K=\Ker \mu$.
    \end{claim}
        Therefore, their $k$-algebra with radical square zero $A(\Lambda)$ and its valued quiver are 
    \begin{equation*}
        A(\Lambda)=l\oplus Z = l \oplus (\Ker (l \otimes_k l \xrightarrowdbl{\mu'} l)), \quad 
        Q_{A(\Lambda)}:
        \begin{tikzcd}
            \bullet \arrow[out=-30,in=30,loop,swap,"{(n-1,n-1)}"]
        \end{tikzcd}
    \end{equation*}
    
    In this example, $\Lambda$ is sg-Hom-finite if and only if $\Lambda$ is Iwanaga-Gorenstein if and only it is Gorenstein if and only if $n=2$.
    \begin{equation}\label{eqa:ex2}
        \begin{tikzcd}[row sep=1.25em, column sep={between origins, 3.5em}]
                &                &                                         &                                           &                                                         &0 \arrow[dd]                                       &                                                      &0 \arrow[dd]                                    &    & \\
                &                &0 \arrow[dd]                             &                                           &0 \arrow[dd]                                             &                                                   &0 \arrow[dd]                                          &                                                &    & \\
                &0 \arrow[rr]    &                                         &K \arrow[rr] \arrow[dl,equal] \arrow[dd]   &                                                         & l \otimes_k xl[\![x]\!]\arrow[rr] \arrow[dl,equal] \arrow[dd]      &                                                      &xl[\![x]\!] \arrow[rr] \arrow[dl,equal] \arrow[dd]   &    &0\\
0 \arrow[rr]    &                &K \arrow[rr,crossing over] \arrow{dd}    &                                           &l\otimes_k xl[\![x]\!] \arrow[rr,crossing over] \arrow[dd]               &                                                   &xl[\![x]\!] \arrow[rr, crossing over] \arrow[dd]           &                                                &0   & \\
                &0 \arrow[rr]    &                                         &K \arrow[rr] \arrow[dl]                    &                                                         &l\otimes_k \Lambda \arrow[rr,"\mu"near start] \arrow[dl] \arrow[dd]             &                                                      &l[\![x]\!] \arrow[rr] \arrow[dl,equal] \arrow[dd] \arrow[dd]   &    &0\\
0 \arrow[rr]    &                &Y \arrow[rr] \arrow[dd]                  &                                           & l\otimes_k l[\![x]\!] \arrow[rr,crossing over] \arrow[dd]             &                                                   & l[\![x]\!] \arrow[rr,crossing over] \arrow[dd]          &                                                &0   & \\
                &                &                                         &                                           &                                                         &l \arrow[rr,"\simeq"near start] \arrow[dl] \arrow[dd]  &                                            &l \arrow[dl,equal]  \arrow[dd]       &    &\\
0 \arrow[rr]    &                &Z \arrow[rr] \arrow[dd]                  &                                           &l \otimes_k l \arrow[rr,"\mu'",crossing over] \arrow[dd]      &                                                   &l \arrow[rr,crossing over] \arrow[dd]  &                                                &0   & \\
                &                &                                         &                                           &                                                         &0                                                  &                                                      &0                                               &    & \\
                &                &0                                        &                                           &0                                                        &                                                   &0                                                     &                                                &    & 
        \arrow[from=2-3, to=4-3, crossing over]  
        \arrow[from=2-5, to=4-5, crossing over]  
        \arrow[from=3-4, to=5-4, "\rotatebox{270}{$\simeq$}" near start]
        \arrow[from=4-3, to=6-3, crossing over]  
        \arrow[from=4-3, to=4-5, crossing over]  
        \arrow[from=4-5, to=6-5, crossing over]  
        \arrow[from=4-7, to=6-7, crossing over]  
        \arrow[from=6-7, to=8-7, crossing over]  
    \end{tikzcd}
    \end{equation}
    ~\\
    \noindent {\textit{Proof of Claim:}}
    Note that $\Gamma \otimes_R \Lambda \simeq l \otimes_k \Lambda$ and $\mu: \Gamma \otimes_R \Lambda \to \Gamma$ induces an isomorphism 
    \begin{equation*}
       \mu: \Gamma \otimes_R \Lambda \otimes_\Lambda (\Lambda/\rad \Lambda) \xrightarrow{\simeq} \Gamma \otimes_\Lambda (\Lambda/\rad\Lambda), 
    \end{equation*}
    where both sides are isomorphic to $l$. So, in this case, $\mu$ is a projective cover and $K$ is a $(\Gamma,\Gamma)$-bimodule. Then by Proposition \ref{prop:sEndGamma_semisimple}, we have 
    \begin{equation*}
        M=\sHom_\Lambda(\Gamma,K) = K/(K\rad \Gamma)
    \end{equation*}
    as a $(D,D)$-bimodule. Moreover, by Proposition \ref{thm:syzygy_in_addGamma} and Corollary \ref{cor:compute syzygy}, we have the big commutative diagram (\ref{eqa:ex2}) and $K=Y\rad \Gamma, Z \simeq Y/(Y\rad \Gamma)$.
    On the other hand, multiplying $x$ induces an isomorphism $Y\rad\Gamma\simeq Y$. Therefore, 
    \begin{equation*}
        M = K/(K\rad\Gamma)= (Y \rad \Gamma / Y (\rad \Gamma)^2)\simeq Y/(Y\rad \Gamma)\simeq Z.
    \end{equation*}
    
\end{example}

\begin{example}\label{ex:3}
    Define an $n\times n$ matrix
    \begin{equation*}
        \Gamma=
        \begin{pNiceMatrix}[last-col]
            R & R & \cdots & R &Q_1 \\
            \pi R & R & \cdots & R &Q_2 \\
            \vdots & \vdots & \ddots & \vdots & \;\, \vdots \\
            \pi R & \pi R & \cdots & R & Q_n 
        \end{pNiceMatrix}
    \end{equation*}
    The labels on the right denote the row vectors which are the indecomposable projective $\Gamma$-modules. 
    
    Consider a partition 
    \begin{equation*}
        J=\Z/n\Z=\{1,2,\dots,n\}=J_1 \sqcup J_2 \sqcup \dots \sqcup J_u.
    \end{equation*}
    
    Define
    \begin{equation*}
        \Lambda=\Lambda[J_1,J_2,\dots,J_u]:=\{(\gamma_{st})\in \Gamma \mid \gamma_{ss}-\gamma_{tt} \in \pi R, \text{ for any $s,t \in J_v$ and $1 \leq v \leq u$}\}.
    \end{equation*}
    Then $(\Lambda,\Gamma)$ is a B{\"a}ckstr{\"o}m order. Define a subset $J'\subseteq J$ by removing all parts of the partition which consist of a single element. This is equivalent to say that $J'=\{j \in J \mid Q_j \not \in \proj \Lambda\}$.

    Let $k_j=Q_j/\rad Q_j\simeq k$ for $j\in J$. 
    For $1\leq v \leq u$, let $P_v$ be the indecomposable projective $\Lambda$-module corresponding to the partition $J_v$.
    The projective cover of a non-$\Lambda$-projective $Q_j$ for $j \in J_v$ is 
    \begin{equation*}
        0 \to \bigoplus_{j'\in J_v \backslash \{j\}} Q_{j'+1} \to P_v \to Q_j \to 0.
    \end{equation*}
    So, we have the semisimple $k$-algebra 
    \begin{equation*}
        D=\sEnd_\Lambda (\Gamma)=\prod_{j\in J'}k_j,
    \end{equation*}
    and as a right $D$-module
    \begin{equation*}
        \sHom_\Lambda(\Gamma, \Omega(Q_j))\simeq \bigoplus_{j'} k_{j'+1},
    \end{equation*}
    where the direct sum is over all $j'\in J_v \backslash \{j\}$ such that $j'+1 \in J'$, or equivalently such that $Q_{j'+1}$ is not a projective $\Lambda$-module. The syzygy functor is given by the diagonal map 
    \begin{equation*}
        \begin{CD}
            \Omega: \quad @. \sEnd_\Lambda (Q_j) @>>> \sEnd_\Lambda (\Omega(Q_j))\\
            @.         @|                           @|                              \\
            \Delta: \quad        @.  k_j                     @>>>  \prod_{j'}k_{j'+1}        
        \end{CD}
    \end{equation*}
    where the product is over all $j'\in J_v \backslash \{j\}$ satisfying the condition as the direct sum above. 
    
    Therefore, the valued quiver $Q_{A(\Lambda)}$ can be obtained combinatorially as follows. For $1 \leq v \leq u$, let $G_v$ be the quiver with vertices $J_v$ and exactly one arrow in each direction between any two distinct vertices. The valued quiver $Q_{A(\Lambda)}$ is obtained from the disjoint union of all $G_v$ for $1 \leq v \leq u$, by 
    \begin{enumerate}[label=\textbullet]
        \item For each arrow from $j \to j'$, change the target $j'$ to $j'+1$. 
        \item Delete the vertices $j$ such that $Q_{j+1}\in \proj \Lambda$ and the arrows adjoining to them.
    \end{enumerate}
    The finite dimensional $k$-algebra with radical square zero $A(\Lambda)$ is the path algebra $kQ_{A(\Lambda)}$ with relations that all paths with length at least two vanish.  
    
    For instance, 
    \begin{enumerate}[leftmargin=16pt]
        \item For $n=2s,s\geq 1$ and $\Lambda=\Lambda[\{1,2\},\{2,3\},\dots,\{2s-1,2s\}]$. Then 
        \begin{equation*}
            Q_{A(\Lambda)}: 
            \begin{tikzcd}[column sep=0.8em]
                1 \arrow[r] & 3 \arrow [r] & \cdots \arrow[r] & 2s-1 & \times & 2 \arrow[out=65,in=110,loop,swap] &  \times & 4 \arrow[out=65,in=110,loop,swap] & \times &  \cdots & \times & 2s \arrow[out=65,in=110,loop,swap]
                \arrow[from=1-4, to=1-1, bend right=45]
            \end{tikzcd}
        \end{equation*}
        and $\Lambda$ is a Gorenstein B{\"a}ckstr{\"o}m order.
        \item For $n=2s+1, s \geq 1$ and $\Lambda=\Lambda[\{1,2\},\{2,3\},\dots,\{2s-1,2s\},\{2s+1\}]$. Then 
        \begin{equation*}
            Q_{A(\Lambda)}: 
            \begin{tikzcd}[column sep=0.8em]
                1 \arrow[r] & 3 \arrow [r] & \cdots \arrow[r] & 2s-1 & \times & 2 \arrow[out=65,in=110,loop,swap] &  \times & 4 \arrow[out=65,in=110,loop,swap] & \times &  \cdots & \times & 2s \arrow[out=65,in=110,loop,swap]
            \end{tikzcd}
        \end{equation*}
        and $\Lambda$ is an Iwanaga-Gorenstein but not a Gorenstein B{\"a}ckstr{\"o}m order.
        \item For $n=2s+1,s\geq 1$ and $\Lambda=\Lambda[\{1,3,\dots,2s+1\},\{2\},\{4\} \cdots,\{2s\}]$. Then 
        \begin{equation*}
            Q_{A(\Lambda)}: 
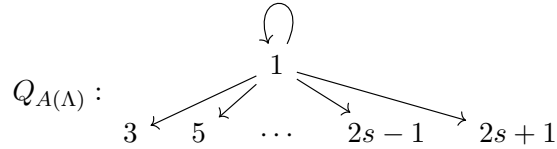
            \begin{tikzcd}[row sep=1em, column sep=1em]
                & & 1 \arrow[out=65,in=110,loop,swap] & & \\
                3 \arrow[urr] & 5\arrow[ur] & \cdots & 2s-3 \arrow[ul] & 2s-1 \arrow[ull] 
            \end{tikzcd} \times 2s+1
        \end{equation*}
        and $\Lambda$ is a sg-Hom-finite but not an Iwanaga-Gorenstein B{\"a}ckstr{\"o}m order.
        \item For $n=2s+2,s\geq 1$ and $\Lambda=\Lambda[\{1,3,\dots,2s+1\},\{2\},\{4\} \cdots,\{2s+2\}]$. Then 
        \begin{equation*}
            Q_{A(\Lambda)}: 
            \begin{tikzcd}[column sep=0.5em]
                1 & \times & 3 & \times & \cdots & \times & 2s+1 
            \end{tikzcd}
        \end{equation*}
        and $\Lambda$ is of global dimension 2.
    \end{enumerate}
\end{example}

\begin{acknowledgement}
    The author would like to thank his supervisor Osamu Iyama for suggesting this problem and for his constant encouragement, insightful discussions, and valuable comments. The author also thanks to Xiao-wu Chen for helpful comments on the paper. The author is also grateful to his labmates, especially Ryu Tomonaga, for their helpful discussions and support throughout his studies.
\end{acknowledgement}

\end{document}